\documentclass{article}
\usepackage{amssymb}

\newtheorem{thm}{Theorem}[section] %

\newtheorem{lemma}[thm]{Lemma} %

\newcommand{\eqn}{\begin{eqnarray}}

\newcommand{\eeqn}{\end{eqnarray}}

\usepackage{amsmath}

\usepackage{graphicx}

\begin{document}

\LARGE{

 \centerline{\textbf{Maslov Indicies in Symplectic Geometry Revisited:}}

\centerline{Sylvain E. Cappell and Edward Y. Miller}

\centerline{Abstract}

  }
  \small{
Thirty years ago in ``On the Maslov Index'' \cite{CLM} the present authors with their late
collaborator Prof. Ronnie Lee of Yale University presented an axiomatization
of the Maslov index, an integer-valued invariant of continuous, piecewise smooth
paths of pairs of Lagrangians, $F: [a,b] \rightarrow (L_1(t),L_2(t))$ in a finite
dimensional symplectic vector space $V$, called here $\mu_{CLM}(\{F(t)\})$. That invariant arises naturally in symplectic geometry and its applications. That paper gave four integer-valued invariants, two geometric and two
analytic, which satisfied these axioms and hence are equal. Here we investigate two
new constructions satisfying those same axioms, so equal to the previous four; they are explicit and more elementary analytic
definitions than those previously considered, each being defined as just an integral
from $a$ to $b$ plus end point corrections.

A variant integer invariant in a different setting, that of a path of symplectic matrices
$M : [a,b] \rightarrow Sp(2n,\mathbb{R})$,
 $\mu_{HZ}( \{ M(t)\})$, was recently defined
by Her and Zhong \cite{HZ}. They described its relation to some previous Maslov indices
in this path context. A main objective of this paper is to clarify the relation between these
two approaches to Maslov indices. In the interest of making this paper self-contained apart
from the previous paper \cite{CLM}, we reiterate definitions of Her and Zhong \cite{HZ} and give an
independent description of the methods of Salomon and Zehnder \cite{SZ1,SZ2}
(which draws on ideas of Gelfand and Lidskii \cite{GL}) which we use.

 Here we also
show that the Her and Zhong invariant may be written as an integral [with a different integrand]
from $a$ to $b$ plus end point corrections. As a result of this, their invariant is here shown to be symplectically
invariant. Given the reference Lagrangian
$L_0=(\mathbb{R}^n \oplus 0)$ in $\mathbb{R}^{2n}$, there is the associated path
of pairs of Lagrangians $t \mapsto  ( L_0, M(t)\cdot L_0)$. In this paper
the difference between $\mu_{HZ}(\{M(t)\})$ and
$\mu_{CLM}( \{t \mapsto  ( L_0, M(t)\cdot L_0) \})$ is explored.
They are here shown to sometimes differ, but in all cases their difference can be written directly
in terms of new end point corrections. This is carried out in several differing
formats. For example, as a simplest case, if the end points $M(a), M(b) \in Sp(2n,\mathbb{R})$
are both in the image of $K: U(n) \subset Sp(2n,\mathbb{R})$, then these two invariants are equal.
At the end of this paper, the relation of these results to the foundational theorems of
Gelfand and Lidskii \cite{GL} is explained.

  }
 \large{

\section{ Statement of Main Results:} \label{sec1}

The title of this paper refers to the following background  . Thirty years ago, in the paper ``On the Maslov index'' \cite{CLM} which we wrote with our late
dear friend and research collaborator Prof. Ronnie Lee of Yale University, we defined
the Maslov indices of symplectic geometry in a variety of settings and showed their equality.
 In particular, it
axiomatized uniquely the Maslov index of a continuous, piecewise smooth path of pairs of Lagrangian subspaces,  $ \mu_{CLM}( \{ \ t \mapsto (L_1(t),L_2(t)) \ \}), \ a\ \le t \le b$,
an integer, in a finite dimensional symplectic vector space, $(V, \ \{ \ , \ \})$
Moreover, we there gave four equivalent definitions, two analytic and two geometric, each
of which satisfied these axioms on paths of pairs of Lagrangians and hence are equal.

\vspace{.1in}
In this paper we offer two additional equivalent definitions of this Maslov index, an integer, here
as an integral from $a$ to $b$ plus a correction at its two end points. They are elementary and  already largely
implicit in \cite{CLM}; the missing ingredient in its proof was lemma \ref{lemma1} below
which extends the unitary context to the general symplectic context.

Also, a parallel formula is given for the recently defined, by Her and Zhong \cite{HZ}, variant Maslov index
of a continuous, piecewise smooth path  $M(t), a  \le t \le b$, in $Sp(2n,\mathbb{R})$,
an integer denoted here by \newline $\mu_{HZ}(\{t \mapsto M(t), a \le t \le b\})$.
The present paper presents formulae for $\mu_{HZ}(\star)$ parallel to the treatment here of
$\mu_{CLM}(\star)$, i.e., as integrals over the path plus end point corrections.

 Given the reference Lagrangian
$L_0=(\mathbb{R}^n \oplus 0)$ in $\mathbb{R}^{2n}$, there is the associated path
of pairs of Lagrangians $t \mapsto ( L_0, M(t)\cdot L_0)$. A goal of this paper is to relate these two integers, $\mu_{HZ}(\{t \mapsto M(t)\})$ and $\mu_{CLM}(\{ t \mapsto (L_0, M(t) \cdot L_0\})$. They are here shown to not be equal in all cases, but in all cases their difference can be written directly
in terms of new end point corrections. This is carried out in several differing
formats. For example, as a simplest case, if the end points $M(a), M(b) \in Sp(2n,\mathbb{R})$
are both in the image of $K: U(n) \subset Sp(2n,\mathbb{R})$, then these two integer invariants are equal.

\vspace{.1in}
In more detail, Her and Zhong in their recent elegant paper ``On the Maslov-type Index for general paths
of symplectic matrices'' \cite{HZ} present a new method of associating
an integer, say $\mu_{HZ}(M)$, to any continuous, piecewise smooth mapping $M:[a,b] \rightarrow Sp(n,\mathbb{R})$
with $a<b$, which serves as a Maslov-type index.
It is related to the  earlier Maslov-type index of Conley-Zehnder \cite{CZ2}, defined
under the conditions $M(a)=  Id(2n)$, the identity matrix, and  $det( M(b) - Id(2n)) \neq 0$. That earlier
Maslov index and its extensions have had extensive applications in the construction of Floer homology \cite{F1,F2}
and the Arnold conjecture \cite{A1,A3, CZ1,CZ2,CZ3,F1,F2,FO,LT,RZ1,RZ2,SZ1,SZ2}.
Because of the importance of Maslov indices, end point modified versions of it and extensions
of its application have been investigated by Liu  \cite{Liu1}, Viterbo \cite{V}, and Long \cite{LO1, LO2} and
Long-Liu, \cite{Liu2} among others; see \cite{HZ} for historical remarks.
All these definitions for paths in the symplectic group have a similar characteristic, the given path is supplemented
by suitable additional paths [three such for the Her and Zhong case with two integrating
to zero] and the integral
taken over the concatenation of these paths being used to define the desired integer invariant.
 Her and Zhong's
paper \cite{HZ} also establishes relations between their generalization, $\mu_{HZ}(\{ t \mapsto M(t)\})$, and some of these earlier versions of this type of Maslov indices. Hence, an integral plus end point corrections

formula, theorem \ref{thmHZ} below,
for the Her and Zhong invariant $\mu_{HZ}(\{ t \mapsto M(t)\})$ may be of interest.

\hspace{.1in}

The following lemma will be helpful in the treatment of this paper.
It addresses a foundational question in symplectic geometry.
Let $Lag(\mathbb{R}^{2n})$ be the space of Lagrangians. As known, the symplectic group
$Sp(2n,\mathbb{R})$ acts transitively on $Lag(\mathbb{R}^{2n})$, say by $M \in Sp(2n,\mathbb{R})$
sends $L$ to $M \cdot L$. Let $K:U(n) \subset Sp(2n,\mathbb{R})$ be the inclusion. It is also
known that the subgroup $U(n)$ also acts transitively on $Lag(\mathbb{R}^{2n})$.  So for a choice of reference Lagrangian
$L_0 \in
Lag(\mathbb{R}^{2n})$ and $M \in Sp(2n,\mathbb{R})$ there is an element $A \in U(n)$
with
$$
K(A) \cdot L_0 = M \cdot L_0.
$$

The question is: What is the simplest relation between $M$ and $A$ when
$M \in Sp(2n,\mathbb{R})$ is written in block form $M= \left( \begin{array}{cc} a & b \\ c & d \end{array} \right) $
with $a,b,c,d$ \ \ $n \times n$ real matrices ?
\begin{lemma} \label{lemma1}
The required result is:
$$ \begin{array}{l}
G[M] = ( a + i \ c) \ (a^t a+ c^t c)^{-1/2} \ \ is \ well \ defined \ and  \ unitary, \ \in U(n),\\
 \ with  \ \ K(G[M]) \cdot L_0= M \cdot L_0.
\end{array}
$$
Also, any  $A \in U(n)$ with $K(A) \cdot L_0 = M \cdot L_0$
is of the form $A = G[M] \cdot w$ with $w \in O(n) \subset U(n)$.
In addition, for $A \in U(n)$ as above,
$$
dim \{ \ kernel \ of\ c \ \} = dim_\mathbb{R} \ L_0 \cap (M \cdot L_0) = dim_\mathbb{R} \ L_0 \cap (K(A) \cdot L_0)= dim_\mathbb{C } \ \{ v \in \mathbb{C}^{2n} \ | \
   (A \cdot A^t ) \cdot v = v \},$$ an equality of geometry and algebra.
\end{lemma}
As seen, $M \mapsto G[M]$ defines a retraction $Sp(2n,\mathbb{R}) \rightarrow U(n)$;
also  the mapping $M \mapsto G[M] \cdot G[M]^t$ induces the well known inclusion, [ see
\cite{CLM}]
$
Lag(\mathbb{R}^{2n}) \subset U(n) \ \ sending \ M \cdot L_0 \ to \ G[M] \cdot G[M]^t.
$

As a direct result of lemma \ref{lemma1}, there is the theorem:
\begin{thm} \label{thmcross}
If $M: [a,b] \rightarrow Sp(2n,\mathbb{R})$ is a continuous, piecewise
smooth mapping, then the equality of Lagrangians
$$
M(t) \cdot L_0 = G[M(t)] \cdot L_0
$$
implies an equality of Maslov indices
$$
\mu_{CLM}(\{ t \mapsto (L_0, M(t) \cdot L_0) \} )=
\mu_{CLM}(\{ t \mapsto (L_0, G[M(t)] \cdot  L_0) \} ).
$$
\end{thm}
This reduction from symplectic to unitary is immediate since one definition of the Maslov index $\mu_{CLM}(\star)$ \cite{CLM}
computes its integer value utilizing only the path of pairs of Lagrangians,
not on how they are expressed.

\vspace{.2in}
Granting lemma \ref{lemma1}, the next theorem appears implicitly in \cite{CLM}.
\begin{thm} \label{thmMain}
Fix the reference Lagrangian $L_0=( \mathbb{R}^n \oplus 0) \subset \mathbb{R}^{2n}$.
Let $M_1, M_2: [a,b] \rightarrow  Sp(2n, \mathbb{R})$ be a pair of continuous, piecewise
smooth mappings.

Let the eigenvalues of $G[M_k(t)] \in U(n), k=1, \ 2$ be written in the form
$$
\{ \lambda_1(t)_k, \cdots, \lambda_n(t)_k\}
$$ where $\lambda_j(t)_k =  \pm \ e^{ i \theta_j(t)_k} $ with $0 \le \theta_j(t)_k < \pi$;
i.e., equivalently $\lambda_j(t)_k^2 = e^{ i (2 \theta_j(t)_k)}$ with $0 \le 2 \theta_j(t)_k < 2 \pi$
for $k=1, \ 2$ respectively.

Let the eigenvalues of $G[M_k(t)] \cdot G[M_k(t)]^t \in U(n)$ be written in the form
$$
\{ \lambda'_1(t)_k, \cdots,\lambda'_n(t)_k\}
$$
where $\lambda'_j(t)_k = e^{i \phi_j(t)_k}$ with $0 \le \phi_j(t)_k < 2 \pi$

Then the Maslov index of the path of pairs of Lagrangians  \cite{CLM},
$$
\mu_{CLM}( \ \{ \ \ t \mapsto (M_1(t) \cdot L_0,  M_2(t) \cdot L_0) \  \})
 = \mu_{CLM}( \{\  t \mapsto (K(G[M_1(t)]) \cdot L_0, \ K(G[M_2(t)]) \cdot L_0) \ \ \} \ )
$$
is equal to
$$ \begin{array}{l}
\int_a^b \ ( t \mapsto det(G[M_2(t)])^2 \in S^1)^\star \ ( \frac{d \theta}{2\pi})
- \int_a^b \ ( t \mapsto det(G[M_1(t)])^2 \in S^1)^\star \ ( \frac{d \theta}{2\pi})\\
 + \Sigma_{j=1}^n \ ( \ \frac{2\theta_j(a)_2 - 2\theta_j(b)_2}{2\pi} \ )
 -  \Sigma_{j=1}^n \ ( \   \frac{2\theta_j(a)_1 - 2\theta_j(b)_1}{2\pi} \ )
 + ( h_2(a)- h_2(b)) - (h_1(a) - h_1(b))\\ \\
 and \ also \ equal \ to \\ \\
 \int_a^b \ ( t \mapsto det(G[M_2(t)])^2 \in S^1)^\star \ (\frac{d \theta}{2\pi})
 - \int_a^b \ ( t \mapsto det(G[M_1(t)])^2 \in S^1)^\star \ (\frac{d \theta}{2\pi}) \\
 + \Sigma_{j=1}^n ( \ \frac{\phi_j(a)_2 - \phi_j(b)_2}{2\pi} \ )
 -  \Sigma_{j=1}^n ( \ \frac{\phi_j(a)_1 - \phi_j(b)_1}{2\pi} \ ) + ( h_2(a)- h_2(b))-( h_1(a)- h_1(b))\\
 \end{array}
 $$
 where $h_k(t)$ is the number of $j$ with $(\lambda_j(t)_k)^2=+1$, or equivalently
 the number of $j$ with $\lambda'_j(t)_k= +1$.

\end{thm}
As seen the Maslov index is expressed in two fashions as an integral
with end point corrections at $a,b$ directly in terms of $M_1(t), \ M_2(t) \in Sp(2n,\mathbb{R})$.
 This remarkable equality
can be traced \cite{CLM} to the equality of the eta invariants
of two distinct self adjoint operators. See section \ref{sec2}.

\vspace{.3in}
 The corresponding
theorem for $\mu_{HZ}( \star) $, the Her and Zhong variant of the Maslov variant \cite{HZ}, is:
\begin{thm} \label{thmHZ}
Let the path $M : [a,b] \rightarrow Sp(2n,\mathbb{R})$ be continuous and piecewise smooth.

Let the $n$ eigenvalues of first kind of Gelfand and Lidskii \cite{GL,SZ1,SZ2} of $M(t)$ be $\lambda\#_1(M(t)), \cdots, \lambda\#_n(M(t))$ with  normalized first kind eigenvalues $(\frac{\lambda\#_j(t)}{|\lambda\#_j(t)|})^2
= e^{ i \theta\#_j(t)}$, with $0 \le \theta\#_j(t) < 2 \pi$.
By definition, the Salomon and Zehnder mapping $\rho: Sp(2n,\mathbb{R}) \rightarrow S^1$ \cite{GL,SZ1,SZ2}
is given by
$$
\rho(M(t)) = \prod_{j=1}^n \frac{\lambda\#_j(t)  }{ | \lambda\#_j(t) | }.
$$
 Let $h\#(t)$ be the number
of $j$ with $\lambda\#_j(M(t))^2 = +1$. Then there
is an equality of integers:
$$ \begin{array}{l}
\mu_{HZ}( \{ t \mapsto M(t) \} )
= \int_a^b \ (t \mapsto \rho(M(t))^2))^\star \ (\frac{d \theta}{2\pi})
  + \Sigma_{j=1}^n \ ( \ \frac{\theta\#_j(a)}{2\pi} - \ \frac{\theta\#_j(b)}{2\pi} \ )
  + (h\#(a)-h\#(b))
  \end{array} $$

  \end{thm}

  In view of the known symplectic invariance of the $n$ eigenvalues of first kind of
an element $M \in Sp(2n,\mathbb{R})$, \cite{SZ1,SZ2}, this theorem implies the strong symplectic
invariance \newline $\mu_{HZ}( \{ t \mapsto M(t) \} ) = \mu_{HZ}( \{ t \mapsto N(t) \cdot M(t) \cdot N(t)^{-1} \} )$
for $N: [a,b] \rightarrow Sp(2n,\mathbb{R})$ continuous and piecewise smooth\footnote{This
invariance should be added to the properties of $\mu_{HZ}(\star)$ listed in \cite{HZ}.}.

A relation between these two theorems is:
\begin{thm} \label{thmsame}
For $ t \mapsto M(t) \in Sp(2n,\mathbb{R}), \ a \le t \le b$ continuous and piecewise smooth
with $M(a) = K(A)$ and $M(b)= K(B)$ for $A, B \in U(n)$, then
$$
\mu_{HZ}( \{t \rightarrow M(t) \}) = \mu_{CLM}( \{t \mapsto (L_0,M(t)   \cdot L_0)\}).
$$
\end{thm}
In section \ref{sec6}, a general formula for the difference
$$
\mu_{HZ}( \{ t \mapsto M(t)\}) - \mu_{V}( \{ t \mapsto (L_0, M(t) \cdot L_0) \})
$$
is given as a sum over the two end points, $a,b$. By example, the  difference is shown not always equal zero.

\vspace{.1in}

For references and insight into the significance of the Maslov index of a path
of pairs of Lagrangians see \cite{CLM, CLM1, CLM2, CLM3}.

\vspace{.3in}
In section \ref{sec2}, background is provided on the integral terms above and details
of the definition of eigenvalues of the first kind are reviewed. In section \ref{sec3}, theorem \ref{thmMain},
an evaluation of $\mu_{CLM}(\star)$, is proved. It consists
mainly in quoting from \cite{CLM}.
The proof of theorem \ref{thmHZ}, an evaluation of the Her and Zhong variant $\mu_{HZ}(\star)$, together with a detailed
exposition of their invariant appears in section \ref{sec4}.
A proof of the algebraic lemma \ref{lemma1} appears in section \ref{secapp}
along with an extended discussion; see theorem \ref{thmalg}.    In section \ref{sec6},
the difference between the distinct invariants, $\mu_{CLM}(\star), \ \ \mu_{HZ}(\star)$  \cite{CLM, HZ}, is made clear
with theorems describing the difference and by an example.
In section 7, it is shown how to recast the integers $\mu_{HZ}(\star), \mu_{CLM}$ in a
much simpler manner. Additionally, the mapping $\Pi : Sp(2n,\mathbb{R}) \rightarrow U(n)$
with $\Pi(K(U)) = U$ for $U \in U(n)$
occurring in the work of Gelfand and Lidskii \cite{GL} also defines an integer by a parallel method.
A variant of this ``Maslov index'' appears in the work of Conley and Zehnder \cite{CZ2}
and also Salmon and Zehnder \cite{SZ2}.
 Sections 8 and 9, gives more details on two foundational theorems of Gelfand and Lidskii \cite{GL}
 and clarifies the difference between using $f : [a,b]\rightarrow Sp(2n,\mathbb{R})$  mapping to  $det_\mathbb{C}(G[f(t)]) \in S^1$ [CLM case] and verses to
 $det_\mathbb{C}(\Pi(f(t))) \in S^1$ [GL case].

\section{Background: The main integral terms:}  \label{sec2}

Since the context of the continuous, piecewise path
$
M : [a,b] \rightarrow Sp(2n,\mathbb{R})
$
directly involves the symplectic group $Sp(2n, \mathbb{R})$,
it is feasible in this paper to just work in the explicit setting of $\mathbb{R}^{2n}$
with its standard symplectic pairing on $2n \times 1$ column vectors $v,w \in \mathbb{R}^{2n}$ defined by
$$
\{ v, w \} = v^t  \cdot J_n \cdot w \in \mathbb{R}.
$$
Here $J_n= \left( \begin{array}{cc} 0 & Id[n] \\ -Id[n]  & 0 \end{array} \right)$
and $Sp(2n, \mathbb{R})$ is the set of $2 n \times 2n$ real matrices $M$ preserving this inner product, i.e.,
 $M^t \ J_n \ M = J_n$. More concretely,
  $M = \left( \begin{array}{cc} a & b \\ c & d  \end{array} \right)$
with $a^t c = c^t a, \ b^t d = d^t b, and \ a^td - c^t b = Id[n]$.

Under the identification $\Phi: \mathbb{R}^{2n} \cong \mathbb{C}^n: \ (x_1,\cdots, x_n, y_1, \cdots, y_n)
\mapsto (x_1+i y_1, \cdots, x_n+i y_n)$, $\{v, w\}$ becomes the imaginary part of the
complex Hermitian inner product $< (z_1, \cdots, z_n),( a_1,\cdots, a_n)>
= \Sigma_{j=1}^n \ \overline{z_j} \cdot a_j$. Under the identification $\Psi$ multiplication
by $i$ on $\mathbb{C}^n$ corresponds to $-J_n$ acting on $\mathbb{R}^{2n}$.

By definition, a Lagrangian subspace $L$ in $\mathbb{R}^{2n}$ is an $n$ dimensional real
subspace for which $\{v,w\}=0$ for $v,w \in L$. Since an element $M$ of $Sp(2n,\mathbb{R})$
preserves $\{\star, \star\}$, multiplication by $M$ induces an bijection to
$Lag(\mathbb{R}^{2n})$, the space of Lagrangians, say $M \in Sp(2n,\mathbb{R})$ sends
$L \mapsto M \cdot L$.

Now for a Lagrangian $L$ we may chose a real basis, say $f_1,\cdots, f_n$, which is additionally
orthonormal for the real symmetric, positive definite inner product $Re   < \star, \star>$. Hence, under the identification
$\Phi : \mathbb{R}^{2n} \cong \mathbb{C}^n$ this basis  become a orthonormal basis of $\mathbb{C}^n$ for the Hermitian
inner product. So, $\Phi(f_1), \cdots, \Phi(f_n)$ are the $n$ column vectors of a element
$A \in U(n)$ and $L$ is the real span of the columns of $A \in U(n)$. Let $L_0$ be the reference Lagrangian $L_0= \mathbb{R}^n \oplus 0$. Alternatively expressed, let $K : U(n) \subset Sp(2n,\mathbb{R})$ be the inclusion
$X + i \ Y \mapsto \left( \begin{array}{cc} X & -Y \\ Y & X \end{array} \right)$,
then $U(n)$ acts transitively on $Lag(\mathbb{R}^{2n})$ via $A \in U(n)$ sends $L_0$
to $K(A) \cdot L_0$ and the induced mapping $U(n) \rightarrow Lag(\mathbb{R}^{2n})$ defines
a surjection
$$
U(n) \rightarrow Lag(\mathbb{R}^{2n}) \ \ which \ \ induces \ a \ bijection \
U(n)/ O(n) \cong Lag(\mathbb{R}^{2n}).
$$
That is, the real span of any orthonormal basis of $\mathbb{C}^n$ under $< \star, \star>$
is a Lagrangian and visa versa.

 Moreover, the mapping
$$ \begin{array}{l}
K(A)  \cdot L_0  \mapsto A \cdot A^t \in U(n)\ \
induces \ an \ inclusion : \ \ Lag(\mathbb{R}^{2n}) \subset U(n).\\
In \ particular, \ the \ mapping \ \
`` \  det(\ )^2 \ "  : Lag(\mathbb{R}) \rightarrow S^1 = \{e^{ i \theta}\}\\
defined \ by \  K(A) \cdot L_0  \mapsto det(A \cdot A^t) = det(A)^2 \in S^1 \ is \ well \ defined.
\end{array}
$$

Also, as $U(n) \rightarrow U(n)/O(n)$ is a fibration, for any continuous,
piecewise smooth mapping $F : [a,b] \rightarrow Lag(\mathbb{R}^{2n})$
there is a continuous, piecewise smooth mapping $A: [a,b] \rightarrow U(n)$
with $F(t)=K( A(t) ) \cdot L_0$.
So for any continuous, piecewise smooth mapping $F: [a,b] \rightarrow Lag(\mathbb{R}^{2n})$
there is the well defined continuous, piecewise smooth mapping
$$ \begin{array}{l}
 t \mapsto `` \ det( \ )^2 \ " \ F(t) \in S^1 \ sending \ \ F(t) = K(A(t)) \cdot L_0 \mapsto det(A(t))^2
 \ for \ A(t) \in U(n). \end{array}
 $$
For example, for any continuous, piecewise smooth mapping $A : [a,b] \rightarrow U(n)$
there is the associated mapping :
$$ \begin{array}{l}
[a,b] \rightarrow S^1 \ \ given \ by  \ t \mapsto det(A(t))^2 \in S^1\\
and \ associated \ integral \ over \ [a,b], \ \
\int_a^b \ \{ \ t \mapsto det(A(t))^2) \ \}^{\star} \ (\frac{d\theta}{2 \pi} ) \in \mathbb{R}
\end{array}
$$

In view of lemma \ref{lemma1} and theorem \ref{thmcross}, for $M: [a,b] \rightarrow Sp(2n,\mathbb{R})$,
the matrix $G[M(t)] \in U(n)$ has the property that
$$
K(G[M(t)]) \cdot L_0 = M(t) \cdot L_0
$$
and hence, since the Maslov index $\mu_{CLM}( \{ L(t) \in Lag(\mathbb{R}^{2n}) \})$ only depends on the path of Lagrangians $L(t), \ a \le t \le b$,
necessarily
$$
\mu_{CLM}( \{ \ t \mapsto M(t) \cdot L_0 \ \}) =  \mu_{CLM}( \{ t \mapsto K(G[M(t)]) \cdot L_0 \}
$$
with respective integral term
$$
\int_a^b \  \{ t \mapsto det(G[M(t)])^2 \}^{\star} \ ( \frac{d \theta}{ 2 \pi} )
$$
appearing as the main constituent of theorem \ref{thmMain}.

\vspace{.3in}

Now we turn to the case of a continuous, piecewise smooth mapping $M : [a,b] \rightarrow Sp(2n,\mathbb{R})$ in the context of the work of Her and Zhong \cite{HZ} summarized below.

For $M \in Sp(2n,\mathbb{R})$, i.e., $M^t \ J_n \ M= J_n$, let $\sigma(M)$ denote
the set of distinct eigenvalues of $M$, say $\lambda_1,\cdots, \lambda_{K(M)}$.
For each $j=1, \cdots, K(M)$, let $Eigen_{\lambda_j}(M)$ be the generalized eigenvalues
of $M$ for this $\lambda_j$, i.e., $Eigen_{\lambda_j}(M) = \{v \in \mathbb{C}^{2n} \ | \
(\lambda_j \cdot Id - M)^{2n}(v)= 0 \}$. Here there is the direct sum decomposition
$ \mathbb{C}^{2n} = \bigoplus_{j=1}^{K(M)} \ Eigen_{\lambda_j}(M)$.

Gelfand and Lidskii \cite{GL} associated to any $M \in Sp(2n,\mathbb{R})$,  the ``$n$
eigenvalues  of the first kind" , say $\{\lambda\#_1(M), \cdots, \lambda\#_n(M)\}$.
 In terms of these, the Salomon and Zehnder \cite{SZ1, SZ2} map $M \mapsto \rho(M)$ is defined by
 by the normalized versions of these first kind eigenvalues
$$
\rho(M) = \prod_{k=1}^n \ \frac{\lambda\#_k(M)}{|\lambda\#_k(M)|} \ \in \ S^1.$$

These first kind eigenvalues are defined following a standard method: Each eigenvalue of the first kind of $M$
is an eigenvalue of $M$. So for each eigenvalue $\lambda_j \in \sigma(M)$ of $M$, the multiplicity
$m_{\lambda_j}(M)$, possibly zero, of the eigenvalues of first kind with eigenvalue $\lambda_j$ is to be defined. These multiplicities are chosen to add up to $n$.

\vspace{.1in}
For $|\lambda_j|<1$ define the first kind multiplicity of $\lambda_j$ as
$ m_{\lambda_j}(M) = dim_\mathbb{C} \ Eigen_{\lambda_j}(M)$. For $|\lambda_j| >1$,
define the first kind multiplicity of $\lambda_j$ as $ m_{\lambda_j}(M) =zero$.

To define the first kind multiplicities for eigenvalues $\lambda_j$ with
$|\lambda_j|=+1$, Gelfand and Lidskii \cite{GL,SZ1,SZ2} [implicitly] introduce
 the complex bilinear pairing
$$ \begin{array}{l}
<<\star, \star>> \ : \ \mathbb{C}^{2n} \times \mathbb{C}^{2n} \rightarrow \mathbb{C} \ \
defined \ by \ << v, \ w>> = (  \ \overline{v}^t \cdot J_n \cdot w \ ).
\end{array}
$$
This bilinear pairing is non-degenerate.

As seen, this inner product is $Sp(2n, \mathbb{R})$ invariant $<< M \ v, M \ w>>= <<v,w>>$
via $M^t \ J_n \ M = J_n$ for $ M \in Sp(2n,\mathbb{R})$. Additionally, if
$M\ v = \lambda \ v$ and $M\ w = \mu \ w$ for eigenvectors $v,w$, then
$$
\overline{\lambda} \mu << v, \ w>>= << M v, \ M w>> = <<v, \ w>>
$$
and so necessarily $<<v, \ w>> =0$ if $ \overline{\lambda} \ \mu \neq +1$. By induction,
this relation extends to the generalized eigenspaces and implies:
If $\lambda, \mu$  are eigenvalues of $M$ with $ \overline{\lambda} \ \mu \neq +1 $,
then $<< Eigen_\lambda(M), \ Eigen_\mu(M) >> = 0$.

Having observed this fact, consider the direct sum decomposition
$$
\mathbb{C}^{2n} = [ \ \bigoplus_{|\lambda_j| \neq 1} \ Eigen_{\lambda_j}\ ] \ \
\bigoplus \ \ [ \ \bigoplus_{|\lambda_k|= 1} \ Eigen_{\lambda_k}\ ]. $$
The bilinear pairing $<<\star , \star>>$ is non-degenerate and, by the above,
the two distinct summands are orthogonal under this pairing; hence, the pairing
restricted to each is non-degenerate on each summand.

Consider the pairing restricted to the second summand,
 $ \bigoplus_{| \lambda_k |= 1} \ Eigen_{\lambda_k}(M)\ ]$ with its induced, non-degenerate pairing.
 But by the above observations, if $\lambda_j= e^{i \theta_j}, \lambda_k= e^{ i \theta_k}$
  with $\theta_j, \  \theta_k$ real and if additionally $\lambda_j \neq
  \lambda_k$, then
  $$
  << \ Eigen_{\lambda_j}(M), \ Eigen_{\lambda_k}(M) \ >> = 0. \ \
$$
In particular, terms of the second direct summand are orthogonal and the induced mapping
for each
$$
<< \star, \ \star>> \ : \ Eigen_{\lambda_j} \times Eigen_{\lambda_j} \rightarrow
  \mathbb{C}$$
  is non-degenerate in this $|\lambda_j|=+1$ context.

Next observe that
  $$
  << iv, \ w>>= - i <<v,\ w>>= - <<v, \ iw>> \ \ (\star)
  $$
and that taking the negative imaginary part  induces the real symmetric [by $(J_n)^t = - J_n$ ] bilinear mapping
  defined by Gelfand and Lidskii and used by Salomon and Zehnder \cite{GL,SZ1,SZ2}:
  $$
  \begin{array}{l}
  \{\{ \star, \ \star\}\} \ : \ Eigen_{\lambda_j} \times Eigen_{\lambda_j} \rightarrow
  \mathbb{R}\\
  \{\{ v, \  w \}\}  = - \ Im \ << v, w>> = - Im \ ( \ (  \overline{v})^t \ J_n \ w \ ).
  \end{array}
  $$
  For these norm $+1$ eigenvalue cases,
    the non-degeneracy of $<< \star, \ \star>>$ coupled to the equalities $(\star)$ implies that the  real symmetric pairing $\{\{ \star, \ \star \}\}$ is also non-degenerate on each
    $ Eigen_{\lambda_j}$.

  The first kind multiplicity of $\lambda_j(M)$ for $|\lambda_j(M)|=+1$,
  namely $m_{\lambda_j(M)}(M)$,
  is now defined to be one half the dimension of the maximal real linear subspace of $Eigen_{\lambda_j}(M) $
  on which $\{\{ v, \  w \}\} $ is positive definite. This maximal real linear
  subspace is of even dimension in view of the equalities $(\star)$.

  The first kind multiplicity of $\lambda_j$ for $|\lambda_j|<+1$, is now defined to
  be $m_{\lambda_j}(M) = dim_\mathbb{C} \ Eigen_{\lambda_j}(M)$;
  The first kind multiplicity of $\lambda_j$ for $|\lambda_j|>+1$, is now defined to
  be zero, $m_{\lambda_j}(M) = 0$.

\vspace{.1in}
  Notable properties of the set of first kind eigenvalues proved by
  Salomon and Zehnder \cite{SZ2} are:

   a) $Sp(2n,\mathbb{R})$ invariance: $M \in Sp(2n,\mathbb{R})$
  and $N \cdot  M \cdot  N^{-1}$ for $N \in Sp(2n,\mathbb{R})$ have the same first kind
  eigenvalues.

   b) If $M= K(A)$ for $A \in U(n)$, then the first kind eigenvalues
  of $M$ are equal to the $n$ eigenvalues of the unitary matrix $A$; all of norm $+1$.

Let there be $L$ distinct eigenvalues of the first kind of $M \in Sp(2n, \mathbb{R})$, say in increasing order:
$$
\lambda_{a_1}, \ \lambda_{a_2}, \cdots, \lambda_{a_L} \ \ with \ 1 \le a_1 < a_3< \cdots <a_L \le 2n.
$$
So by the above algorithm each $\lambda_{a_t}$ has multiplicity $m_{\lambda_t}(M)$
eigenvalues of the first kind with the value $\lambda_{a_t}$ for $t=1,\cdots, L$.

   They proceed
  to give an axiomatic treatment of the associated mapping
  $$
  \rho : Sp(2n,\mathbb{R}) \rightarrow S^1 \ \ defined \ by \
  M \mapsto   \prod_{t=1}^{L}  ( \frac{\lambda_{a_t}}{|\lambda_{a_t}|} \ )
  ^{ m_{\lambda_t}(M)   } \in S^1. $$

 Hence, there is an associated
integral for a continuous, piecewise smooth $M: [a,b] \rightarrow Sp(2n,\mathbb{R})$
given by
$$ \begin{array}{l}
\int_a^b \  ( t \mapsto  \rho(M(t))^2)^\star \  (\frac{d \theta}{  2\pi}).
\end{array}
$$
This is the main constituent of $\mu_{HZ}( \{t \mapsto M(t)\} )$.
Note the square.

 \section{Proof of theorem \ref{thmMain}:} \label{sec3}
 Theorem \ref{thmMain} is implicit in the work \cite{CLM}.

Given $M_1, M_2 : [a,b] \rightarrow S(2n,\mathbb{R})$, then by theorem \ref{thmcross}
there is an equality
$$
\mu_{CLM}(\{ t \mapsto ( M_1(t) \cdot L_0, \  M_2(t) \cdot L_0) \})
=  \mu_{CLM}(\{ t \mapsto ( K(G[M_1(t)])  \cdot L_0, \ K(G[ M_2(t)]) \cdot L_0) \}).
$$
This equality allows us to reduce the computation of the Maslov index of
a path of pairs of Lagrangians to the special case of
$$
\mu_{CLM}(\{ t \mapsto ( K(G[M_1(t)])  \cdot L_0, \ K(G[ M_2(t)]) \cdot L_0) \}).
$$
Then the methods of \cite{CLM} which assumed such a unitary reduction apply
can be cited without loss of generality; replacing $M_j(t) \in Sp(2n,\mathbb{R})
$ by $G[M_j(t)] \in U(n)$.

Theorem \ref{thmMain} above, may be compared to Theorem 0.4 on page 124 of \cite{CLM}
which states in the present terminology:
\begin{thm}[Theorem 0.4]  \label{thmCLM} For $H(t) = ( L_1(t), L_2(t)), \  a \le t \le b$, \ \
a continuous, piecewise path of pairs of Lagrangians with $L_k(t) = K(A_k(t)) \cdot L_0$
for continuous, piecewise smooth mappings $t \mapsto A_k(t) \in U(n) \ with \  k=1,2$,
$$ \begin{array}{l}
\mu_{CLM}( t \mapsto  (L_1(t), \ L_2(t))) = \int_a^b [ L_2^\star(\omega) - L_1^\star(\omega)] \\
+ \frac{1}{2} \ ( \ \eta( D(L_1(b),L_2(b)))  - \eta(D(L_1(a), L_2(a))) \ )
+ \frac{1}{2} ( h(a) - h(b)).
\end{array}
$$
with $\eta(D(\star,\star))$ the respective eta invariant with specified boundary conditions
and $h(t) = dim_\mathbb{R} (L_1(t) \cap L_2(t))$.
\end{thm}
Here $\int_a^b \ L_k^\star (\omega)$ is a short hand for
$\int_a^b \  ( t \mapsto  det(A_k(t))^2 \ )^\star \ ( \frac{d\theta}{2 \pi} )  \in \mathbb{R}$ with $k=1,2$.

As explained in \cite{CLM}, the eta invariant may be taken as either of two distinct
self adjoint operators with equal eta invariants.

Let $L_1,L_2$ be Lagrangian subspaces
of $\mathbb{R}^{2n}$ with $L_2 = K(A) \cdot L_1$ for $A \in U(n)$. Define a real operator on the space of smooth vector valued function
$\Phi: [0,1] \rightarrow \mathbb{R}^{2n}$ with boundary conditions
$$
\phi(0) \in L_1 \ and \ \phi(1)\in L_2,
$$
It is given by $D(L_1,L_2) \ \phi = -J \frac{\phi(t)}{dt}$. This gives a self adjoint operator
$D(L_1,L_2)$ with kernel identified with $L_1 \cap L_2$. Let $A$ have eigenvalues
$\lambda_1(A), \cdots, \lambda_n(A)$ which are written in the form
$$
\lambda_j(A)= \pm \cdot e^{i \theta_j(A)}   \ with  \  0 \le \theta_j(A) < \pi.
$$
In \cite{CLM} the associated eta invariant $\eta(D(L_1,L_2))$ is shown to be equal to
$$
\eta(D(L_1,L_2)) = \Sigma_{ 0 < \theta_j(A) <\pi} \ \ ( \  1 - 2 ( \frac{\theta_j}{\pi}) \ )
$$

Secondly, define the real self adjoint operator $D^\#(L_1,L_2) =-i \ \frac{d\psi}{dt}$
acting on the smooth mappings $\psi : [0,1] \rightarrow \mathbb{C}^{2n}$ with boundary\
condition
$$
\psi(1)= ( A \cdot A^t ) \cdot \psi(0).
$$

In \cite{CLM} it is shown that $kernel(D\#(L_1,l_2)) =  \{ v \in \mathbb{C}^{2n} \ | \
(A \cdot A^t) \ v = v\} = ( L_1 \cap L_2) \otimes \mathbb{C}$.
 In particular,
the real dimension of the kernel of $D(L_1,L_2)$ equals the complex dimension
of the kernel of $D\#(L_1,L_2)$. Let the eigenvalues of $A \cdot A^t$ be
$\lambda\#_1,\cdots, \lambda\#_n$ with $\lambda\#_j(A)= e^{i \phi_j(A)}$ with $0 \le
\phi_j(A) < 2 \pi$. Moreover, the eta invariant $\eta(D\#(L_1,L_2))$
is shown to equal
$$
\eta(D\#(L_1,L_2)) = \Sigma_{0 < \phi_j< 2\pi} \ (1 - 2 ( \frac{\phi_j}{2 \pi})
$$

Inserting these results with $A_k(t) = G[ M_k(t)], k=1,2$ into theorem \ref{thmMain}
gives its proof.

\section{Proof of theorem \ref{thmHZ}:} \label{sec4}

Proof of theorem \ref{thmHZ}:

Let the eigenvalues of the first kind for $M(t)$
be $\{\lambda\#_1(t), \cdots, \lambda\#_n(t)\}$ with associated normalized eigenvalues of
the first kind $\{ \frac{\lambda\#_j(t)}{ | \lambda\#_j(A)|  } \}$, where
$ \frac{\lambda\#_j(t)}{ | \lambda\#_j(t)|}= e^{ i \phi\#_j(t)}$ are given uniquely
by conventions to be discussed later.
 Let
$Diag( t)$ denote the diagonal block
sum of the associated $2 \times 2 $ rotation matrices
$$
C(\phi_j(t)) = \left(   \begin{array}{cc}  cos(\phi\#_j(t)) & - sin( \phi\#_j(t)) \\
                         sin(\phi\#_j(t)) & cos(\phi\#_j(t))  \end{array} \right)=
                         K[ (e^{i \phi\#_j(t)}) ],  \ j=1, \cdots, n,
                         $$
                         namely, iterating: the diagonal block sum is:
                         $$ Diag(t) = C(\phi\#_1(t)) \diamond  C(\phi\#_2(t)) \diamond \cdots \diamond C(\phi\#_n(t))=
                         K( (e^{i \phi\#_1(t)}) \diamond  (e^{i \phi\#_2(t)}) \diamond \cdots \diamond (e^{i \phi\#_n(t)})) \ \in \ K(U(n)).$$

                            With,  for   example, diagonal block sum:
 $$ \left( \begin{array}{cc} a & b \\ c & d \end{array} \right)
                         \diamond \left( \begin{array}{cc} A & B \\ C & D \end{array} \right)
                        \stackrel{def.}{ =}  \left( \begin{array}{cccc} a & 0 & b & 0\\
                          0 & A & 0 & B\\
                          c & 0 & d & 0 \\
                            0 & C& 0 &D \end{array} \right).
                            $$

                         In particular, the first kind eigenvalues of this sum of diagonal blocks
                         consist of  exactly the complex norm one eigenvalues
                         $\{ e^{i \phi\#_j(t)}, j=1, \cdots n \}$ and so by definition are also the normalized first kind eigenvalues
                         of $M(t)$.

Following Her and Zhong \cite{HZ},
 consider the end $t=a$ and the block diagonal $2n \times 2n$ matrix,
 $Diag(a)$
in some fixed order. By $Sp(2n,\mathbb{R})$ path connected,
we may chose a smooth path, say
$$
P[a]:  [ a-1,a] \rightarrow Sp(2n,\mathbb{R}) \ \ with \ \ P(a-1)=Diag(a),  \ \ P(a) = M(a)
$$
from $Diag(a) $ to $M(a)$, Since the normalized eigenvalues of these end points
are equal, the mapping $\rho(\star)$ carries them to the same point of $S^1$.
Hence, the associated mapping $S^1 \rightarrow S^1$ has some degree, say $D(a)$.
Now it is known, see \cite{CZ1,CZ2}, that the loop defined by
$$
Diag( e^{2 \pi  J_n t}, 1, \cdots, 1), \ \ 0 \le t \le 1
$$
maps by $\rho(\star)$ to $S^1$ with degree $+1$; hence introduce the loop
$$ \begin{array}{l}
{}L[a] : [a-2, a-1] \rightarrow Sp(2n,\mathbb{R}) \ via \ \
 {} L[a]((a-2)+t)  = Diag(e^{- 2 \pi  D(a) J_n  t},1, \cdots, 1) \cdot Diag(a) \end{array},
$$
then the concatenation $L[a] \star P[a] : [a-2, a] \rightarrow Sp(2n,\mathbb{R})$
composed with $\rho$
has degree zero. In particular, $$\int_{a-2}^a  ( t \mapsto \rho(L[a]\star P)(t))^2)^\star
( \frac{d \theta}{2\pi})  = 2\int_{a-2}^a  ( t \mapsto \rho(L[a]\star P)(t))^\star
( \frac{d\theta}{2\pi}) =0.$$ By these means, a smooth path from
$Diag(a)$ to $M(a)$, $L[a] \star P[a]$, has been defined with integral  of the pull back via $\rho(\star)^2$ equal to zero.

 Let this smooth path from $Diag(a)$ to $M(a)$ be denoted by $Initial :[a-2, a]
\rightarrow Sp(2n,\mathbb{R})$ with $Initial(a-2)= Diag(a), Initial(a)= M(a)$
with $\int_{a-2}^{a} \ ( t \rightarrow \rho(Initial(t))^2)^\star ( \frac{d \theta}{2\pi} )=0$.

In a similar manner we may construct a smooth path $Final:[b,b+2]$ with
$Final(b)= M(b), \newline Final(b+2) = Diag(b)$ and $ \int_{b}^{b+2} \ ( t \rightarrow \rho(Final(t))^2)^\star ( \frac{d \theta}{2\pi} ) =0$.

By this method, the path $M$ may be extended to the concatenation
$Total = Initial \star M \star Final : [a-2,b+2] \rightarrow Sp(2n,\mathbb{R})$
with $Total(a-2)= Diag(a), Final(b+2)= Diag(b)$ with \newline
$\int_{a-2}^{b+2} \ ( t \rightarrow \rho(Total(t))^2)^\star ( \frac{1}{2\pi} ) ( d \theta)=
\int_{a}^{b} \ ( t \rightarrow \rho(M(t))^2)^\star ( \frac{ d \theta}{2\pi} ) $. This
replaces the ends $M(a),M(b)$ by the block diagonal sums $Diag(a), Diag(b)$
which just record the normalized eigenvalues of the first kind
without any change of integral.

Now there is a difficulty in having the end points of our path with eigenvalues
$\pm 1$, so multiply the above path $M(t)$  apparatus by $e^{+ J_n \epsilon}$
getting $t \mapsto e^{+J_n \epsilon} \cdot Total(t)$.  Now consider the end points
of $Total$ under this change:
By the nature of the diagonal block sum above, [recall multiplication by $i$ on $\mathbb{C}^n$
corresponds to applying $-J_n$ on $\mathbb{R}^{2n}$]
$$
e^{ J_n \epsilon} \cdot Diag(t) = e^{ J_n \epsilon} \cdot K( (e^{i \phi\#_1(t)}) \diamond  (e^{i \phi\#_2(t)}) \diamond \cdots \diamond (e^{i \phi\#_n(t)}))= K( (e^{i (\phi\#_1(t)- \epsilon)}) \diamond  (e^{i (\phi\#_2(t)- \epsilon)}) \diamond \cdots \diamond (e^{i (\phi\#_n(t) - \epsilon)})).$$

Thus for $\epsilon >0$ sufficiently small,
both of $ e^{+J_n \epsilon} \cdot Diag(a)$  and   $e^{J_n \epsilon} \cdot Diag(b)$ will have
eigenvalues all distinct from $\pm 1$. Here the explicit diagonal block orthogonal nature
of these matrices is used to see the effect of multiplication by $e^{+J_n \epsilon}$.

Now it is necessary to relate in some manner  the block diagonal end point matrices $ e^{J_n \epsilon} \cdot Diag(b)$ to $e^{J_n \epsilon} \cdot Diag(a)$
by a path of some sort.
Actually, Her and Zhong connect  $e^{J_n \epsilon} \cdot Diag(b)$ to a diagonal block
matrix $e^{J_n \epsilon} \cdot Diag\#(a)$ with $(e^{J_n \epsilon} \cdot Diag\#(a))^2 =(e^{J_n \epsilon} \cdot Diag(a))^2$. Their method is described next.

 For this
it is useful to uniquely parameterize
the normalized eigenvalues of the first kind of $M(b),M(a)$ in the following \textbf{unconventional} manner:  [Remark: This differs from that used in theorem \ref{thmHZ}, so at the end
we must convert to the more standard representation used in theorem \ref{thmHZ}.]
$$
\begin{array}{l} \frac{\lambda\#_j(a)}{|\lambda\#_j(a)|} = e^{ i (A_j(a)+ \pi \cdot r_j(a)} \ with \
    0 < A_j(a) \le \pi, \ r_j(a) \in \{0,1\} \\
    \frac{\lambda\#_j(b)}{|\lambda\#_j(b)|} = e^{ i (A_j(b)+ \pi \cdot r_j(b)} \ with \
    0 < A_j(b) \le \pi, \ r_j(b) \in \{0,1\} \\ \end{array} $$
 Note the restrictions on the angles are the half open intervals,
 $  (0,\pi ], \ (\pi , 2 \pi ]$. For example, $+1$ is recorded
 as $+1 = e^{2 \pi i} =e^{i ( \pi + \pi \cdot 1)} $ with parameters $(\pi,1)$
 and $-1$ is recorded as $-1 = e^{i \pi} = e^{i( \pi + \pi \cdot 0)}$ with parameters $(\pi,0)$.
 $e^{i \phi}$ with $0 < \phi < \pi$ is recorded by $(\phi,0)$ and
 $e^{i(\phi+ \pi)}$ with $0 < \phi <\pi$ is recorded by $(\phi,1)$.

 Now suppose $\epsilon>0$ is  chosen small so that the only $
 \frac{\lambda\#_j(a)}{|\lambda\#_j(b)| }, \frac{\lambda\#_j(b)}{|\lambda\#_j(b)| }$ except those equal to $+1$ or $-1$ are within  the  angular interval
  $2 \epsilon$ of these values $\pm 1$. Then $e^{-i\epsilon} \lambda_j(a)$ and $e^{-i \epsilon} \lambda_j(b)$
  have all their parameters
  $$
  (r,s) \ \ with \  \ r\in  [ \epsilon, \pi -\epsilon] \sqcup [ \pi + \epsilon, 2\pi - \epsilon]
  \ \ and \ \ s \ equal \ to \ 0, \ 1 \ respectively.
  $$
  Here $s$ equals $r_j(a) \in \{0,1\}$ for $t=a$ and equals  $r_j(b)  \in \{0,1\}$ for $t=b$.

  [The idea is that those with values $+1$ are pushed down by $\epsilon$ to have angle $- \epsilon= 2 \pi - \epsilon$
  and those with values $-1$ are pushed down to have angles $\pi-\epsilon$;
  and the other whose angles are in the range $2 \epsilon $ to $\pi - 2 \epsilon$
  or range $\pi+ 2 \epsilon$ to $2\pi - 2 \epsilon$ will be pushed down by $\epsilon$
  to have range  $ \epsilon$ to $\pi -  \epsilon$
  or range $\pi+ \epsilon$ to $2\pi- \epsilon$.]

At this point, Her and Zhong \cite{HZ} appeal to a lemma proved by Salmon and Zehnder \cite{SZ1} and utilized
by Long \cite{LO1,LO2}:

\begin{lemma} Let $Sp_L(2n, \mathbb{R})$ denote the set of elements of
$Sp(2n,\mathbb{R})$ having $L \in \mathbb{C}$ as an eigenvalue.
Then the set
$$
X = Sp(2n,\mathbb{R}) - Sp_1(2n,\mathbb{R})-Sp_{-1}(2n,\mathbb{R})
$$
of elements without eigenvalues $\pm 1$ is a union of finitely
many connected path components each of which are simply
connected.
\end{lemma}

They add to the points
$$ \begin{array}{l}
Diag(b) = C(  e^{ i (A_1(b) + \pi \cdot r_1(b))} ) \diamond C( e^{ i (A_2(b) + \pi \cdot r_2(b))} )
\diamond \cdots \diamond C( e^{ i (A_n(b) + \pi \cdot r_n(b))}  ),  \\
Diag(a) = C( e^{ i (A_1(a) + \pi \cdot r_1(a))} ) \diamond C( e^{ i (A_2(a) + \pi \cdot r_2(a))} )
\diamond \cdots \diamond C( e^{ i (A_n(a) + \pi \cdot r_n(a))} ),  \\
the \ additional \ point \ of \ Sp(2n,\mathbb{R}), \\
Diag(a)\# = C( e^{ i (A_1(a) + \pi \cdot r_1(b))} ) \diamond C( e^{ i (A_2(a) + \pi \cdot r_2(b))} )
\diamond \cdots \diamond C( e^{ i (A_n(a) + \pi \cdot r_n(b))} \ ).  \\
\end{array}
$$

\vspace{.1in}
They observe that the two diagonal block matrices $ e^{J_n \epsilon} \cdot Diag(b), e^{J_n \epsilon} \cdot Diag(a)\#$ are in $X$ and moreover $ e^{J_n \epsilon} \cdot Diag(a),  e^{J_n \epsilon} \cdot Diag(a)\#$ lie in the same connected component of $X$.
Hence, they may define a smooth path, say $P : [b+2, b+3] \rightarrow Sp(2n,\mathbb{R}) $
in $X$ with $P(b+2)=  e^{J_n \epsilon} \cdot Diag(b)$ and $P(b+3) =e^{J_n \epsilon} \cdot Diag(a)\#$. As the components of $X$ are
simply connected,  any two such smooth paths with the same end points will have equal
integrals. Also, the path is homotopic relative to its end points under the added
constraint that it lies in $X$.

Having done so, the concatenation $H: total \star P : [a-2, b+3] \rightarrow Sp(2n, \mathbb{R})$
defines a continuous, piecewise smooth mapping with $H(b+3)^2 = ( e^{J_n \epsilon} \cdot Diag(a)\#)^2 = (e^{J_n \epsilon} \cdot Diag(a))^2=H(a-2)^2$.
Hence, the mapping $t \mapsto \rho(H(t))^2$ results in a mapping of the closed loop
mapping to $S^1$:
$Loop := [a-2,b+3]/ (a-2) \sim (b+3) \stackrel{ t \mapsto \rho(H(t))^2}{\rightarrow } S^1$.
[As $\rho(Y)^2 = \rho(Y^2)$ for the $Y$'s of diagonal orthogonal block type.] The degree of
the mapping $S^1 \cong Loop \rightarrow S^1$ defined by $t \mapsto  \rho( H(t))^2$
is defined to be their variant \textbf{Maslov index} $\mu_{HZ}(G)$.
The degree is evaluated by an integral, namely :

$$
\begin{array}{l} \int_{a-2}^{b+3} ( t \mapsto \rho(H(t))^2)^\star (\frac{d \theta}{\pi}) ( d\theta)
  \ \ \ \  \ =  degree \ of \ mapping \ of \ circles \\\stackrel{def.}{=} \mu_{HZ}(\{t \mapsto H(t), \  a \le t \le b\}).
  \end{array}
  $$

\vspace{.3in}
Now consider the observation from the above, that the angles $A_j(b) + \pi \cdot r_j(b)- \epsilon $ for $e^{ J_n \epsilon} \cdot Diag(b)$
and $A_j(a) + \pi \cdot r_j(b) - \epsilon $  for $ e^{J_n \epsilon} \cdot Diag(a)\#$ lie in the \textbf{same} closed interval
$$
\pi \cdot r_j(b) + \epsilon \ \le \ A_j(b) + \pi \cdot r_j(b), \ A_j(a) + \pi \cdot r_j(b) \
\le  \ (\pi \cdot r_j(b) +\pi) - \epsilon
$$
for each $j$. So they are easily connected by the homotopy deforming
the angles via \newline $ (1-t) \cdot (A_j(b) + \pi \cdot r_j(b) - \epsilon ) + t (\cdot A_j(a) + \pi \cdot r_j(b) - \epsilon)\ for \
\ 0 \le t \le 1$
inside the convex closed interval $[ \ ( \pi \cdot r_j(b) - \epsilon), \ (\pi \cdot  r_j(b) +\pi - \epsilon) \ ] $ for $j=1, \cdots, n$. In particular, none of the eigenvalues of the associated
diagonal block matrices  are equal to $\pm 1$.
Hence, the thus  associated special  path, say $P(t)$, is in diagonal block form and  lies in $X$ and shows that $ e^{ J_n \epsilon} \cdot Diag(b) $ and $ e^{ J_n \epsilon} \cdot Diag(a)\#$ are in the same component as desired to use in the Her and Zhong definition of $\mu_{HZ}( \star)$.

But for this special choice  $P(t)$ through $X$, the required integral from
$ e^{ J_n \epsilon} \cdot Diag(b)$ to $e^{ J_n \epsilon} \cdot Diag\#(a)$ is easily
computed as:
$$ \begin{array}{l}
\int_{b+2}^{b+3} ( t \mapsto \rho(P(t))^2)^\star \ ( \frac{d\theta}{2\pi})
 = \Sigma_{j=1}^n  \frac{ (A_j(a) + \pi \cdot r_j(b) -\epsilon) -( A_j(b) + \pi \cdot r_j(b)
 - \epsilon)}{\pi}\\

 =  \Sigma_{j=1}^n  \frac{ 2A_j(a)}{2 \pi} -  \Sigma_{j=1}^n \frac{2A_j(b)}{2 \pi}\\
 \end{array}
 $$

Now compare the conventions on angles in theorem \ref{thmHZ}, namely $\theta\#_j(t)$
with the choices above $\phi\#_j(t)$.
Here,  in the second instance, $\frac{\lambda\#_j(t)}{ |\lambda\#_j(t)|}= e^{ i(A_j(t) + \pi r_j(t))}$ with $0 < A_j(t) \le \pi$
and $r_j(t) \in \{0,1\}$; so $(\frac{\lambda\#_j(t)}{|\lambda\#_j(t)|})^2 = e^{i(2 A_j(t))} $ with $0 < 2A_j(t) \le  2 \pi $
to be compared with $(\frac{\lambda\#_j(t)}{|\lambda\#_j(t)|})^2 = e^{ i \theta\#_j(t)}$ with $0 \le \theta\#_j(t) < 2$.
 Hence, if $(\frac{\lambda\#_j(t)}{|\lambda\#_j(t)|})^2 \neq +1$, then
$2 A_j(t) = \theta\#_j(t)$; and if  $(\frac{\lambda\#_j(t)}{|\lambda\#_j(t)|})^2 =+1$, then
$2 A_j(t) =2\pi$ while  $\theta\#_j(t)=0$. Thus making the translation into
the terms of theorem \ref{thmHZ}
$$\begin{array}{l}
\int_{b+2}^{b+3} ( t \mapsto \rho(P(t)))^\star \ ( \frac{d\theta}{\pi})
=  \Sigma_{j=1}^n  \frac{ 2A_j(a)}{2 \pi} -  \Sigma_{j=1}^n \frac{2A_j(b)}{2 \pi}\\
= \Sigma_{j=1}^n \frac{ \theta\#_j(a) - \theta\#(b)}{2\pi} + (h\#(a) - h\#(b))
\end{array}
$$
Therefore:
$$ \begin{array}{l}
\mu_{HZ}( \ \{ t \mapsto M(t), a \le t \le b \} \ ) =
\int_{a-2}^{b+3} ( t \mapsto \rho( P \star Total)(t))^2)^\star  ( \frac{d\theta}{2 \pi})\\
= 0 + \int_a^b \ (t \rightarrow \rho(M(t))^2)^\star   ( \frac{d\theta}{2 \pi})
 + 0 +  \Sigma_{j=1}^n \frac{ \theta\#_j(a) - \theta\#(b)}{2\pi} + (h\#(a) - h\#(b))
\end{array}
$$
as claimed.

\section{Proof of Lemma \ref{lemma1} and theorem \ref{thmalg}: } \label{secapp}
Let $M \in Sp(n,\mathbb{R})$ be symplectic, say
$
M = \left( \begin{array}{cc} a & b \\ c & d \end{array} \right);
$
i.e.,   $ M^t \cdot J_n \cdot M = J_n$ or more explicitly
$
a^t c = c^t a  \ and \  b^td = d^t b  \ and \ a^t d - c^t b = Id[n].
$

\vspace{.2in}
Now define $D[M] = a^t \ a + c^t \ c$. It is real symmetric, positive definite,
so we may define $D[M]^{-1/2}$.

Reason: For $v \in \mathbb{R}^{2n} \neq 0$ an eigenvector with real eigenvalue $\lambda$.
$$
\begin{array}{l}
\lambda \ v^t v = v^t \cdot D[M] \cdot v = v^t \ a^t \ a \ v + v^t \ c^t \ c \ v
= (av)^t \ (av) + (cv)^t \ (cv) = |av|^2 + |cv|^2 \ge 0
\end{array}
$$
which is real and non-negative. But if $\lambda =0$, then $av=0, \ cv=0$
which imply that $v^t a^t =0, v^t c^t=0$, which implies
that $v^t \ v = v^t Id[n] v  = v^t \ ( a^t d - c^t b) \ v = 0$. So $\lambda=0$
is not an eigenvalue.

\vspace{.1in}
\begin{thm} \label{thmalg}
Let $K :U(n) \subset Sp(2n,\mathbb{R})$ be the inclusion $X+ i \ Y \mapsto
\left( \begin{array}{cc} X & -Y \\ Y & X \end{array} \right)$. Here $X^t Y = Y^tX$ and  $X^tX+Y^tY=Id[n], \ i \mapsto -J_n$.

\begin{itemize}
\item[1.] For $M \in Sp(2n,\mathbb{R})$ as above; the $n \times n$ \textbf{complex matrix}
$$
G[M] = ( a + i \ c) \cdot D[M]^{-1/2}
$$
is unitary.

Note that for $(X+i \ Y) \in U(n)$, $G[ K(X+i \ Y)]=(X + i \ Y)$. Hence, $M \mapsto G[M]$ is a retraction of
$Sp(2n,\mathbb{R})$ onto $U(n)$.

\item[2.] Let $L_0 = \mathbb{R}^n \oplus 0$ be the reference Lagrangian and
$M \cdot L_0$ be the result of the action of $M$ on $L_0$. That is,
$M \cdot L_0$ is the real span of the first $n$ columns of $M$.
Then
$$
M \cdot L_0 = K(G[M]) \cdot L_0.
$$
That is, the span of the first $n$ columns of $M$, $=\left( \begin{array}{c|c} a  & b \\ c &d \end{array} \right)$,
 equals the span of the first $n$ columns of $K(G[M])$,
$ =\left( \begin{array}{c|c} a \cdot D[M]^{-1/2} & - c \cdot D[M]^{-1/2}\\ c \cdot D[M]^{-1/2} &  a \cdot D[M]^{-1/2} \end{array} \right)$.

For example, for any continuous path $M: [a,b] \rightarrow Sp(2n,\mathbb{R})$
there is the equality
$$
\{M(t) \cdot L_0 = K(G[M(t)]) \cdot L_0, \ a \le t \le b\}:
$$
a reduction from symplectic to unitary geometry.

\item[3.] The symplectic matrix $M $ equals a product of
of symplectic matrices
$$
M = K(G[M]) \cdot \left( \begin{array}{cc} D[M]^{+1/2} & 0 \\ 0 & D[M]^{-1/2} \end{array} \right)
\cdot \left( \begin{array}{cc} Id[n] & \alpha \\ 0 & Id[n] \end{array} \right)
$$
with $\alpha = D[M]^{-1} \cdot (a^t b + c^td)$. Moreover, $\alpha$ is symmetric; i.e.,
$$
 D[M]^{-1} \cdot (a^t b + c^td)  =( b^ta +d^t c) \cdot D[M]^{-1}.
 $$

 \item[4.] There is a natural deformation retraction of $Sp(2n,\mathbb{R})$ onto $U(n)$
 defined for $M \in Sp(2n,\mathbb{R})$ and $0 \le u \le 1$ by :
 $$
 H(M,u) =K(G[M]) \cdot \left( \begin{array}{cc} D[M]^{+1/2} & 0 \\ 0 & D[M]^{-1/2} \end{array} \right)^{(1-u)}
\cdot \left( \begin{array}{cc} Id[n] & \alpha \\ 0 & Id[n] \end{array} \right)^{(1-u)}
 $$

\end{itemize}
\end{thm}

\vspace{.3in}
Proof: Using $a^t c = c^t a$:
$$ \begin{array}{l}
(\overline{G[M]})^t \cdot G[M] = D[M]^{-1/2} \ (a^t - i \ c^t) \ (a+  i \ c) \cdot D[M]^{-1/2}\\
 =  D[M]^{-1/2}\ (\ (a^ta + c^t c) + i ( a^tc- c^ta) \ ) \cdot D[M]^{-1/2}\\
 =  D[M]^{-1/2} \  D[M] \ D[M]^{-1/2} = Id[n]\\
 \end{array} $$

Item 2 is clear since $D[M]^{-1/2}$ defines an isomorphism $\mathbb{R}^n \rightarrow
\mathbb{R}^{n}$.

Item 3 is a direct calculation: with $\alpha =D[M]^{-1} \cdot (a^t b + c^td)$. Here
$G[M] \in U(n)$ so $G[M]^{-1} = \overline{G[M]}^t = ( (a \cdot D[M]^{-1/2})^t - i \
(c \cdot D[M]^{-1/2})^t) = D[M]^{-1/2}(a - i\  c)^t $. Here $-i \mapsto J_n
= \left( \begin{array}{cc} 0 & 1 \\ -1 & 0 \end{array} \right)$.
 Hence,
$$
\begin{array}{l}
K(G[M])^{-1} \cdot M  = \left( \begin{array}{cc} D[M]^{-1/2} & 0 \\ 0 & D[M]^{-1/2} \end{array}
\right) \cdot \left( \begin{array}{cc} a^t & c^t \\ -c^t & a^t \end{array} \right) \cdot
\left( \begin{array}{cc} a & b \\ c & d \end{array} \right) \\
 \left( \begin{array}{cc} D[M]^{-1/2} & 0 \\ 0 & D[M]^{-1/2} \end{array} \right) \cdot \left( \begin{array}{cc} a^ta+ c^tc & a^tb+ c^td\\
   -c^ta+a^tc & -c^tb+a^td \end{array} \right) \\
   =\left
   ( \begin{array}{cc} D[M]^{-1/2} & 0 \\ 0 & D[M]^{-1/2} \end{array} \right) \cdot \left( \begin{array}{cc} D[M] & a^tb+ c^td \\ 0 & Id[n] \end{array} \right)  \\ =
    \left( \begin{array}{cc}  D[M]^{+1/2} & D[M]^{-1/2} ( a^tb+ c^td )\\ 0 & D[M]^{-1/2} \end{array} \right)\\

    = \left( \begin{array}{cc} D[M]^{+1/2} & 0 \\ 0 & D[M]^{-1/2} \end{array} \right) \cdot

       \left( \begin{array}{cc} Id[n] & D[M]^{-1} ( a^tb+ c^td ) \\ 0 & Id[n] \end{array} \right)

                                       \end{array} $$
Since the matrices $M, K(G[M])$, and $\left( \begin{array}{cc} D[M]^{+1/2} & 0 \\ 0 & D[M]^{-1/2} \end{array} \right)$ are symplectic, necessarily $ \left( \begin{array}{cc} Id[n] & \alpha \\ 0 & Id[n] \end{array} \right)$ is simplectic. Hence, $\alpha^t = \alpha$.

This completes the proof of item 3. Item 4 is clear as the powers are well defined.

\section{Comparison of $\mu_{HZ}( \star)$ and $\mu_{CLM}( \star)$; An example of
their divergence:}  \label{sec6}

The equality appearing in theorem \ref{thmsame}
$$
\mu_{HZ}( \{t \rightarrow M(t) \}) = \mu_{CLM}( \{t \mapsto (L_0, M(t) \cdot L_0)\})
$$
for $M: [a,b] \rightarrow Sp(2n,\mathbb{R})$ with $M(a)= K(A), \ M(b)= K(B)$
with $A,B \in U(n)$
is proved as follows:

Theorem \ref{thmalg} part 4, defines a continuous, piecewise smooth, strong deformation retraction:
 $$H: Sp(2n, \mathbb{R}) \times [0,1] \rightarrow Sp(2n, \mathbb{R})$$
of $Sp(2n, \mathbb{R})$ onto $U(n)$. That is, $H(M,0)= M$, $H(K(A),u)= A$ for $ A \in U(n)$ and $0 \le u \le 1$,
and $H(M, 1) = G[M] \in U(n)$.

Hence, there is the associated continuous, piecewise smooth mapping
$$ \begin{array}{l}
S : [a,b] \times [0,1] \stackrel{H(M(t),u)}{\rightarrow} Sp(2n,\mathbb{R}) \stackrel{\rho(\star)^2}{\rightarrow} S^1 \
defined  \ \ by \ \
S(t,u) = \rho(H(M(t),u))^2 \in S^1. \end{array}.$$

Now restricted to the bottom edge $u=0$, the pull back of $( \frac{d \theta}{2 \pi})$ integrated
gives the integral part of \newline $\mu_{HZ}(\{ t \mapsto M(t) \} )$ appearing in theorem
\ref{thmHZ}. When the end point corrections are added, it equals this Maslov index
$\mu_{HZ}(\star)$.

 While restricted to the top edge $u=1$, the pull back of $( \frac{d \theta}{2 \pi})$
integrated
gives the integral part of $\mu_{HZ}(\{ t \mapsto K(G[M(t)]) \} )$. But $G[M(t)]\in U(n)$
so $\rho(K(G[M(t)])^2 = det(G[M(t)])^2$ so this second integral is exactly that appearing
as the integral part of \newline $\mu_{CLM}(\{ t \rightarrow ( L_0, K(G[M(t)])\cdot L_0)\})$
which by theorem \ref{thmcross} equals
$ \mu_{CLM}(\{ t \rightarrow ( L_0, M(t)\cdot L_0)\})$. When the end point corrections are
added, it equals this Maslov index, $\mu_{CLM}(\star)$.

But by $M(a)= K(A), M(b)=K(B)$ for $A, B \in U(n)$, the paths obtained by
setting $t=a$ are the constant paths and for $t=b$ are the constant paths; so contribute
zero to the integral over the boundary of the rectangle $[a,b] \times [0,1]$. Also, the end point corrections cancel here.

But since the map $S$ extends from the boundary to the whole rectangle, the
closed 1-form $( \frac{d \theta}{2 \pi})$ integrates to zero along the whole boundary
and putting in the eight end point corrections, the sum is again zero.

Consequently, this shows that the integral for $\mu_{HZ}(\{ t \mapsto M(t) \} )$
and for $\mu_{CLM}(\{ t \rightarrow ( L_0, M(t)\cdot L_0)\})$ are equal. But
since the end points are in the image of $U(n)$ the normalized first kind eigenvalues
for $K(A)$ equal the eigenvalues of $A$, and similarly for $K(B)$ and $B$,
 the end point corrections for each of these two paths are also equal, so
$$
\mu_{HZ}(\{ t \mapsto M(t) \} ) = \mu_{HZ}(\{ t \mapsto (L_0, M(t) \cdot L_0 ) \} )
$$
follows. This proves theorem \ref{thmsame}.

\vspace{.1in}
The same method applies to any continuous, piecewise smooth function
$
M : [a,b] \rightarrow Sp(2n,\mathbb{R})$
giving  the result:
\begin{thm} \label{thmgen1}
With $H(t,u)$ defined as in theorem \ref{thmsame},
$$ \begin{array}{l}
\mu_{HZ}( \{ t \mapsto M(t) \} ) - \mu_{CLM}( \{ t \mapsto (L_0, M(t) \cdot L_0 \} ) \\
= \mu_{HZ}( \{ t \mapsto M(t) \} ) - \mu_{CLM}( \{ t \mapsto (L_0,K( G[M(t)])  \cdot L_0 \} ) \\
= \mu_{HZ}( \{ u \mapsto H(a, u) \in Sp(2n,\mathbb{R}),\ u=0,1  \})
 - \mu_{HZ}( \{ u \mapsto H(b, u) \in Sp(2n,\mathbb{R}),\ u=0,1  \}). \end{array}
   $$
   \end{thm}
The difference only depends on the end points $M(a), M(b) \in Sp(2n,\mathbb{R})$.
The integrals and their corrections for the sides $t=a$ and $t=b$ contribute
the resulting corrections.

\vspace{.3in}
Another formula of this type relies on a result of Conley and Zehnder, section 1.5, page 218,
`` \ Normal forms for distinct eigenvalues \ " \cite{CZ2}. It states: ``Assume $M \in Sp(2n,\mathbb{R})$
has distinct eigenvalues, none of which is equal to $-1$. Then, as is easily verified, there is a symplectic basis, in which the matrix has block diagonal form. Every block corresponds to an eigenvalue group
and has one of the following three normal forms, where we use the abbreviation
$$
R(\alpha) = \left( \begin{array}{cc}  cos(\alpha) & - sin(\alpha) \\
 sin(\alpha) & cos(\alpha) \end{array} \right).
 $$
 \begin{itemize}
 \item[1.] Hyperbolic plane (eigenvalue group $(\beta, \beta^{-1})$, \  $\beta \neq 0$ real.
 $
 M[1] = \left( \begin{array}{cc} \beta & 0 \\ 0 & \beta^{-1} \end{array} \right)$.
 Here the normalize eigenvalue of the first kind is $sign(\beta)$ and $G[M[1]]$
 is the $1 \times 1$ matrix $sign(\beta)$.

 \item[2.] Ellipitic plane $(\lambda, \overline{\lambda}, \ \lambda = e^{i \alpha}$, $\alpha $ real
 and $\alpha \neq \pm 1$.
 $ M[2] = R(\alpha)$. Here the normalized eigenvalue of the first kind is $e^{i \alpha}$
 and $G[M[2]] $ is the $1 \times 1$ matrix $e^{i \alpha}$.

\item[3.] Complex eigenvalue group $( \lambda, \lambda^{-1}, \overline{\lambda}, (\overline{\lambda})^{-1}$, \ \ $\lambda= \rho \cdot e^{i \theta}$, with $|\rho| <1$
    and $\lambda \neq \pm 1$.
    $M[3]= \left( \begin{array}{cc} \rho \ R(\theta) & 0 \\ 0 & \rho^{-1} \ R(\theta)
    \end{array} \right)$. Here the normalized eigenvalues of the first kind are
    $\{  e^{i \theta},  e^{-i \theta}\}$ while the eigenvalues of $G[M[3]]$ are
    $\{  e^{i \theta}, e^{i \theta} \}$\ ".
    \end{itemize}

As observed in each case the normalized eigenvalues of $M[k]$ are equal to the eigenvalues
of the unitary $G[M[k]]$. Hence, as conjugation preserves the eigenvalues of a unitary
matrix and preserves the normalized eigenvalues of the first kind for a symplectic matrix,
any $M\in Sp(2n, \mathbb{R})$ with distinct eigenvalues and with no eigenvalue equal
to $-1$, has its normalized eigenvalues of the first kind equal to the eigenvalues
of the unitary form $G[M]$.

By this result, for any $M \in Sp(2n,\mathbb{R})$ with distinct eigenvalues and
no eigenvalue equal to $-1$, there is an element, say $X \in Sp(2n,\mathbb{R})$
with
$$
X_M \cdot M \cdot X_M^{-1} = diagonal\ block \ form \ sum \ of \ such \ 3 \ types
$$

So pick a smooth path $P: [0,1] \rightarrow Sp(2n,\mathbb{R})$ with $P(0)= Id[2n]$
and $P(1)= X_M$ and apply the above approach. Getting
a mapping
$$ \begin{array}{l}
T: [0,1] \times [0,1] \rightarrow Sp(2n, \mathbb{R}) \\
T(t,u)=\rho( H( P(t) \cdot M \cdot P(t)^{-1}, u))^2 \in S^1
\end{array}
$$

Now restricting to the bottom edge, $u=0$ we have the constant
path $\rho( P(t) \cdot M \cdot P(t)^{-1})^2 = \rho(M)^2$ since the
normalized first kind eigenvalues are invariant under conjugation. Also
the end point corrections cancel.

Now restricting to the top edge $u=1$, gives the result
$\rho( K(G[P(t) \cdot M \cdot P(t)^{-1}]))^2 = det( G[P(t) \cdot M \cdot P(t)^{-1}])^2$
so is the integral part of $\{t \mapsto (L_0,G[P(t) \cdot M \cdot P(t)^{-1}] \cdot L_0)\}
=  \{ t \mapsto (L_0, (P(t) \cdot M \cdot P(t)^{-1}) \cdot L_0) \}$.
That is, the integral part of $\mu_{CLM}(\{ t \mapsto (L_0, (P(t) \cdot M \cdot P(t)^{-1}) \cdot L_0) \} ).$ Restoring the end point corrections gives the Maslov index.

Now restricting to the edge $t=1$ the deformation is a constant one as the result
from first kind of eigenvalues are constant under the deformation of this
sum of diagonal blocks of only three types. So the
integral from this piece is zero and the  end point corrections vanish. There is
no contribution.

The remaining edge $t=0$ reclaims the value for
$  \{ u \mapsto \rho( H(M, u))^2 \}$ above.
In combination, adding in the non-canceling
end point corrections we have proved:

\begin{thm} \label{thmlast}
 Assume $M \in Sp(2n,\mathbb{R})$
has distinct eigenvalues, none of which is equal to $-1$.
Pick as ia possible a smooth path $P: [0,1] \rightarrow Sp(2n,\mathbb{R})$
with $P(0)=Id[2n]$
with
$$
P(1) \cdot M \cdot P(1)^{-1}
$$
 a diagonal block sum of the above three standard types.

Define the mapping
$$ \begin{array}{l}
W: [0,1] \times [0,1] \rightarrow Sp(2n, \mathbb{R})\\
by  \ W(t,u) = H( P(t) \cdot M \cdot P(t)^{-1}, u ) \in Sp(2n,\mathbb{R})
\end{array}
$$
Then [via two sides cancel integrals and end point correction,
the remaining integrals and associated end point corrections add to zero :
$$
\mu_{CLM}(\{ t \mapsto (L_0, (P(t) \cdot M \cdot P(t)^{-1}) \cdot L_0) \}
=\mu_{HZ}( \{ u \mapsto  H( M,u)\} )
$$

\end{thm}

Applied to theorem \ref{thmgen1}, this yields:

\begin{thm} \label{thmgen2}
For $M : [a,b] \rightarrow Sp(2n, \mathbb{R})$
such that $M(a)$ and $M(b)$
have distinct eigenvalues and none equal to $-1$,
then as possible, there exists exist two continuous,
piecewise smooth maps
$$
X, \ Y : [0,1] \rightarrow Sp(2n,\mathbb{R})\
with \ \ X(0)= Id[2n], \ Y(0)= id[2n]
$$
such that $X(1) \cdot M(a) \cdot X(1)^{-1}$ and
$Y(1) \cdot M(b) \cdot Y(1)^{-1}$ are a diagonal
block sum of the three above standard types.

For such a choice,
$$
\begin{array}{l}
\mu_{HZ}( \{ t \mapsto M(t)\} ) - \mu_{CLM}( \{ t \mapsto (L_0, M(t) \cdot L_0\} )\\
= \mu_{CLM}( \{ u \mapsto ( L_0, (X(u) \cdot M(a) \cdot X(u)^{-1}) \cdot L_0)\})
 -  \mu_{CLM}( \{ u \mapsto ( L_0, (Y(u) \cdot M(b) \cdot Y(u)^{-1}) \cdot L_0)\})
 \end{array}
 $$
 \end{thm}

\vspace{.1in}
An example of the deviation of the two indices, $\mu_{HZ}( \{ t \mapsto M(t) \in Sp(2n,\mathbb{R}), \
a \le t \le b \} )$ and $\mu_{V}( \{ t \mapsto (L_0, \ M(t)\cdot L_0) \
a \le t\le b\} )$ is easily obtained.

Let $$
Rot(\theta)=\left( \begin{array}{cc} cos(\theta) & -sin(\theta) \\ sin(\theta) & cos(\theta)
\end{array} \right) \in K(U(2)) \\
and \ M(S) = \left( \begin{array}{cc} S & 0 \\ 0 & S^{-1} \end{array} \right) \in Sp(2n, \mathbb{R})
$$
for $\theta$ real varying and $S>1$ real fixed.

Define $F: [-\pi/4, + \pi/4] \rightarrow Sp(2, \mathbb{R})$ by $F(\theta) = R(\theta) \cdot M(S) \cdot
R(\theta)^{-1} \in Sp(2,\mathbb{R})$. Then the Her and Zhong index
$$
\mu_{HZ} ( \{ t \rightarrow F(t), \ -\pi/4 \le t \le +\pi/4\}) \ vanishes
$$
as by symplecitic invariance $\rho(F(\theta))=\rho(M(R)) \in Sp(2,\mathbb{R}$ is constant.

While :
$$
\begin{array}{l}

 F(\theta)
     =  \left( \begin{array}{cc} S \ cos^2(\theta) + S^{-1} \ sin^2(\theta) & cos(\theta) sin(\theta)
      (S - S^{-1}) \\ cos(\theta) sin(\theta)  (S-S^{-1}) & S \ sin^2(\theta) + S^{-1} \ cos^2(\theta) \end{array} \right) \end{array}
     $$
 Hence, $F(\theta) \cdot L_0 = \{ T \cdot \left( \begin{array}{l} S  cos^2(\theta) + S^{-1}  sin^2(\theta) \\  cos(\theta) sin(\theta)  (S-S^{-1}) \end{array} \right), \ T \ real \}$.
 So for $-\pi/4 \le \theta \le + \pi/4$, \ \  $ L_0 \cap  (F(\theta) \cdot L_0) = L_0$ only if
  $cos(\theta) sin(\theta) =0$, i.e., $\theta =0$ where $F(0) \cdot L_0 = L_0$.
  But the linear approximation to the first column of $F(\theta) $ at $\theta=0$ is just
  $ \left( \begin{array}{l} S \\ (S-S^{-1}) \Delta \theta \end{array} \right)$;
  so the Lagrangian $F(\theta) \cdot L_0$ crosses $L_0$ for $\theta=0$ transversally
  in the positive direction. Thus, $\mu_{CLM}( f : [-\pi/4, +\pi/4] \rightarrow Sp(2,\mathbb{R})
  =+1$ while $\mu_{HZ}( f : [-\pi/4, +\pi/4] \rightarrow Sp(2,\mathbb{R}) =0$.

\section{ Simpler Formulae :} \label{sect7}

Given a continuous mapping $f: [a,b] \rightarrow S^1 = \{ e^{i \theta}\}$,
there is a natural integer $N(f)$ `counting'' the signed number of times $f(t)$ crosses
the point $1 = e^{i \cdot 0}$ which generalizes the degree of the induced map
$[a,b]/ \ a\sim b \rightarrow S^1$
for the case that $f(a)=f(b)$. It employs a convention specified below.

\vspace{.1in}
Also, in the bel this simple method will be shown to be relevant to the Maslov indices,
$\mu_{HZ}( \star), \mu_{CLM}(\star)$,that of Conley and Zehnder \cite{CZ2},
and also Salomon and Zehnder \cite{SZ2}.

\vspace{.2in}
Given such continuous $f : [a,b] \rightarrow S^1$, let $\Psi: \mathbb{R} \rightarrow S^1$ be the universal covering mapping
$$
\Psi(\theta)= e^{i \theta} \in S^1.
$$
Thus, having chosen an element $\theta_{f}(a) \in \mathbb{R}$ with $\Psi( \theta_f(a)) = f(a)$,
there is a unique continuous covering mapping $\theta_f: [a,b] \rightarrow \mathbb{R}$ with
$$
\Psi( \theta_f(t)) = f(t) \in S^1 \ \ for \ t \in [a,b].
$$

The homotopy
$$  \begin{array}{l}
H: [a,b] \times [0,1] \rightarrow \mathbb{R} \ \ defined \ by \\
H(t,u) =(1-u)\cdot ( (1-(\frac{t-a}{b-a}  )) \ \theta_f(a) + (\frac{t-a}{b-a}) \ \theta_f(b))
 + u \cdot \theta_f(t)
 \end{array}
 $$
 continuously deforms the path $\theta_f(t)$ in $\mathbb{R}$ for $u=1$ to the constant speed path
 $ ( (1-(\frac{t-a}{b-a})) \ \theta_f(a) + (\frac{t-a}{b-a}) \ \theta_f(b))$ for $u=0$ leaving
 the end points $H(a,u)= \theta_f(a),\ H(b,u)= \theta_f(b)$ fixed for all $u$.

 Hence, it suffices
 to describe the crossing number at $+1 \in S^1$ of the image of the  constant speed path  $( (1-(\frac{t-a}{b-a})) \ \theta_f(a) + (\frac{t-a}{b-a}) \ \theta_f(b))$ for $t$ ranging
 from $0$ to $1$ passing from $\theta_{f(a)}$ to $\theta_{f(b)}$ in $\mathbb{R}$.

 If $\theta_f(a)) = \theta_f(b)$, the integer $N(f)$ is declared to be zero.
 If $ \theta_f(a) < \theta_f(b)$, the integer $N(f)$ is declared to be the
 number of integers $K$ for which $\frac{\theta_f(a)}{2\pi} \le K  <  \frac{\theta_f(b)}{2 \pi}$.
 If  $ \theta_f(a) > \theta_f(b)$, the integer $N(f)$ is declared to be minus  the
 number of integers $K$ for which $\frac{\theta_f(b)}{2\pi} \le K  <  \frac{\theta_f(a)}{2 \pi}$.

 [Note the special convention of end points.]

 As seen if $f(a)= f(b)$, then $N(f)$ equals the required degree. If $f: [a,b] \rightarrow S^1$
 and $g: [b, c] \rightarrow S^1$ have $f(b)=g(b)$, then the concatenation
 $(f\star g) : [ a,c] \rightarrow S^1$ with $(f \star g) | [a,b] = f$
 and $(f\star g)|[b,c] = g$ has $N(f\star g)= N(f) + N(g)$. Moreover,
 by the adopted convention
 if $\epsilon>0$ is chosen sufficiently small, then for any $\eta$ with $0 < \eta \le \epsilon $
 there is the equality
 $$
 N(f) = N( e^{-i \eta} \ f) \in \mathbb{Z}.
 $$.

If $ f : [a,b] \rightarrow S^1$ is piecewise smooth, then
$$
 \frac{\theta_f(b)}{2 \pi} - \frac{\theta_f(a)}{2 \pi} = \int_a^b  \ ( t \mapsto f(t) \in S^1)^\star \ \frac{d \theta}{2 \pi}.
 $$
 In the continuous case, it is still quite natural to define
 $$
  \int_a^b  \ ( t \mapsto f(t) \in S^1)^\star \ \frac{d \theta}{2 \pi}
  $$
  by $\frac{\theta_f(b)}{2 \pi} - \frac{\theta_f(a)}{2 \pi}$.

 Similarly, if $f$ is piecewise smooth,
 and  $\theta\#(t)$ be uniquely defined by $f(t) = e^{i \theta\#(t)}$ with
$0 \le \theta\#(t) < 2 \pi$ and $h\#(t)=1$ if $f(t)=1$ and zero otherwise. Then
$$
N(f) =
 \int_a^b \ (t \mapsto f(t) \in S^1))^\star \ (\frac{d \theta}{2\pi})
  +  ( \ \frac{\theta\#(a)}{2\pi} - \ \frac{\theta\#(b)}{2\pi} \ )
  + (h\#(a)-h\#(b)) \in \mathbb{Z}.   \ \ (\star \star)
  $$
  Also with the above convention about the integral, for $f$ continuous
this equality still holds.

\vspace{.3in}
More generally if $f: [a,b] \rightarrow (S^1)^n$ is continuous, then we may define the integer
$N(f)$ by
$$
N(f) = \Sigma_{j=1}^n \  N( \pi_j(f)(t) \ \ with \ \ f(t)= (\pi_1(f)(t), \cdots, \pi_n(f)(t).
$$

 \vspace{.1in}
 In these terms, $\mu_{CLM}( t \mapsto (L_0, M(t) \cdot L_0))$ may be re-expressed as follows:
 Consider the induced continuous mapping $ t \mapsto G(M(t)) \in U(n)$ and its associated eigenvalues
 as a continuous function of $t$, all of norm $+1$, say $\lambda_j(G(M(t)), j=1,\cdots, n$. These may be arranged in a continuously
 ordered fashion and for each $j$ a continuous covering function $\theta_j(t) : [a,b] \rightarrow
 \mathbb{R}$ chosen with
 $$
 e^{i \theta_j(t)} = \lambda_j(G[M(t)]) \in S^1 \ \ for \ t \in [a,b].
 $$
Then if $f$ is piecewise smooth, there is the equality
$$
\mu_{CLM}( t \mapsto (L_0, M(t) \cdot L_0)) = \Sigma_{j=1}^n \ N(t \mapsto \lambda_j(G[M(t)])^2)
\in \mathbb{Z}
$$
based on $det_C( G(M(t))^2$ and the conventions of \cite{CLM}.
This is just a restatement for $M_1(t)=Id, \ M_2(t)=M(t)$ of the equality of theorem \ref{thmMain}.

\vspace{.3in}
Similarly for $\mu_{HZ}(M : [a,b] \rightarrow Sp(2n,\mathbb{R}))$, one has
for each $t \in [a,b]$ the $n$, norm $+1$, normalized first kind eigenvalues of $M(t)$,
say $\lambda_j(M(t))/|\lambda_j(M(t))|$, varying continuously.
Then if $f$ is piecewise smooth, there is the equality
$$
\mu_{HZ}( t \mapsto  M(t)) = \Sigma_{j=1}^n \ N(t \mapsto (\lambda_j(M(t)/|\lambda_j(M(t))|)^2)
\in \mathbb{Z}
$$
Again Her and Zhong's conventions, adapted from those of \cite{CLM}, are those
appearing here and in theorem \ref{thmHZ}.

\vspace{.1in}
There is a analogous mapping
$$
\Pi: Sp(2n, \mathbb{R})  \rightarrow U(n) \ with \ \ K(\Pi(M)) = M \cdot (M^t \ M)^{-1/2}
$$
defined by Gelfand and Lidskii \cite{GL}. For continuous
 mapping $M : [a,b] \rightarrow Sp(2n,\mathbb{R})$, there is defined via the eigenvalues
 of $\Pi(M(t)$ a Maslov type index. A variant of this appears in the work of
 Conley and Zehnder \cite{CZ1} and Salomon and Zehnder \cite{SZ2}.

\section{Added Relations to the work of Gelfand and Lidskii \cite{GL}} \label{sect8}

From the last section, the integrals appearing in the Maslov index of \cite{CZ1,SZ2}
are of the type
$$
\int_a^b \ ( f \rightarrow  det_C(\Pi(M(t)))^\star \ \frac{d \theta}{2\pi}
$$
while those of \cite{CLM} are of the form
$$
\int_a^b \ ( f \rightarrow  det_C(G[M(t)])^2)^\star \ \frac{d \theta}{2\pi}
= 2 \ \int_a^b \ ( f \rightarrow  det_C(G[M(t)]))^\star \ \frac{d \theta}{2\pi}
$$
for the retractions $\Pi , G : Sp(2n,\mathbb{R}) \rightarrow U(n)$.

In theorem \ref{thmalg} above, there is an identification for $M(t)=
\left( \begin{array}{cc} a(t) & b(t) \\ c(t) & d(t) \end{array} \right) \in Sp(2n, \mathbb{R})$
as a product of $K[M(t)] \in U(n)$ with and a supplementary term:
$$
M(t)= K(G[M(t)]) \cdot \left( \begin{array}{cc} D[M(t)]^{+1/2} & 0 \\ 0 & D[M(t)]^{-1/2} \end{array} \right)
\cdot \left( \begin{array}{cc} Id[n] & \alpha(t) \\ 0 & Id[n] \end{array} \right)
$$
with $\alpha(t) = D[M(t)]^{-1} \cdot (a(t)^t b(t) + c(t)^t(d))$. Moreover, $\alpha(t)$ is symmetric.

In the work of Gelfand and Lidskii \cite{GL} there is a parallel theorem for $\Pi: Sp(2n,\mathbb{R})
\rightarrow U(n)$, namely:
 Let $J = \left( \begin{array}{cc} 0 & 1 \\ -1 & 0 \end{array} \right)$.

A foundational theorem of Gelfand and Lidskii \cite{GL} is:
\begin{thm} \cite{GL} \label{thmGL1}
For $N$ a real $m \times m$ matrix with $det(N) \neq 0$.  Let
$$
N = P(N) \cdot Q(N)
$$
be the unique polar decomposition of $N$. That is, $Q(N) = \sqrt{N^t \cdot N}$, a real symmetric positive
definite matrix, and $P(N) = N \cdot Q(N)^{-1}$ which is a real orthogonal matrix, i.e.,
$P(N)^t  \cdot P(N) = Id[m]$.

\vspace{.1in} If $M(t) \in Sp(2n,\mathbb{R})$, then
$P(M(t))$ and $Q(M(t))$ are also symplectic. Indeed, for $K: U(n) \subset Sp(2n, \mathbb{R})$
the inclusion, $P(M(t)) = K( \Pi(M(t)))$ for unique $\Pi(M(t)) \in U(n)$ and $Q(M(t))$ is uniquely
expressed in the form
$$
Q(M(t)) = e^{A(M(t))} \ \ with \ \ A(M(t)) = \left( \begin{array}{cc} a_1(M(t)) & a_2(M(t)) \\ a_2(M(t)) & -a_1(M(t)) \end{array}
\right)
$$
for real symmetric $n \times n$ matrices $a_1(M(t)), \ a_2(M(t))$. Here \newline
 $ J \   \left( \begin{array}{cc} a_1(M(t)) & a_2(M(t)) \\ a_2(M(t)) & -a_1(M(t)) \end{array} \right)
=  - \ \left( \begin{array}{cc} a_1(M(t)) & a_2(M(t)) \\ a_2(M(t)) & -a_1(M(t)) \end{array} \right) \cdot J$
implies that $e^{A(M(t))} \in Sp(2n,\mathbb{R})$.

Additionally, under the standard identification \newline $C: \{n \times n \ real \  symmetric\} \cong \mathbb{R}^{n(n-1)/2}$
$$
M(t) \mapsto (\Pi(M(t)),  C(a_1(M(t))), C(a_2(M(t)))) \in U(n) \times  \mathbb{R}^{n(n-1)/2} \times \mathbb{R}^{n(n-1)/2}
$$
is a smooth bijection.
\end{thm}

Combining these theorems, we see that precise relation of the retractions $M(t) \rightarrow K[M(t)] \in U(n)$ and $M(t) \rightarrow \Pi(M(t)) \in U(n)$. In theorem \ref{thmalg} $M(t)$ is written
in the form $K[M(t)] \cdot (\star)$ while in theorem \ref{thmGL1} $M(t)$ is written in the
form $\Pi(M(t)) \cdot (\star)$.

In the next section, in combination of another foundational theorem of Gelfand and Lidskii \cite{GL},
the relation of the integrands
$$
( f \rightarrow  det_C(\Pi(M(t))^2)^\star \ \frac{d \theta}{2\pi} \ \ and \ \
( f \rightarrow  det_C(G[M(t)]^2)^\star \ \frac{d \theta}{2\pi}
$$
will be computed.

\section{Another Foundational theorem of Gelfand and Lidskii \cite{GL}
and its application}   \label{sect9}

Consider the theorem:

 \begin{thm} \cite{GL} \label{thmGL2} Let $f: [a,b] \rightarrow Sp(2n,\mathbb{R})$ be a smooth
 mapping, then the mapping $H: [a,b] \rightarrow H(t)$ defined by
 $$
 H(t) = -J \ \frac{ d \ f(t)}{dt} \cdot f(t)^{-1} \ \ for \ \ a \le t \le b
 $$
 is a smooth mapping to the real symmetric $2n \times 2n$ matrices.
 Here $f(t)$ satisfies the first order differential equation
 $$
 \frac{f(t)}{dt} = J \cdot H(t) \cdot F(t).   \ with \ H(t) \ real \ symmetric. \ (\star)
 $$
 Conversely, if $H(t)$ is real symmetric and smoothly varying,
 then the unique solution to $(\star)$ satisfying the initial
 condition $f(0)= Id[2n]$ has solution with $f(t) \in Sp(2n,\mathbb{R})$
 for all $a \le t \le b$.
 \end{thm}
[ \  Actually they also considered an extended version replacing $Sp(2n,\mathbb{R})$
 by the group  \newline $\{M \ | \ M \ 2n \times 2n \ complex\ with\ \
 (\overline{M})^t \ J \ M = J \}$ and show that
 for smooth maps, say $f\#$,  to this group $-J \ \frac{ d \ f\#(t)}{dt} \cdot f\#(t)^{-1}$
 is Hermitian. \ ]

 In the theorem's context, since the eigenvalues of the real symmetric matrix $H(t)$ are real, we get a remarkable mapping
 $$
 t \mapsto Trace( -H(t)) = Trace(  J \ (\frac{ d \ f(t)}{dt}) \cdot f(t)^{-1}) \in \mathbb{R}
 $$
 recording the sum (all real) of the eigenvalues of $-H(t)$.

 The first theorem addresses unitary items:
\begin{thm} \label{thm4}
  Let $g: [a,b] \rightarrow U(n)$ be a smooth mapping to the unitary matrices
  which is smoothly diagonalizable by unitary matrices.
 Say,
 for some smooth mapping $V :[a,b] \rightarrow U(n)$
 $$
 g(t)= V(t) \cdot Diag(t) \cdot V(t)^{-1}
 $$
 where $Diag(t)$ is a  diagonal matrix with $Diag_{j,j}(t) = e^{i \theta_j(t)}, \ j=1,\cdots, n$
 for smooth real valued functions $\theta_j(t), \ j=1,\cdots, n$ defined for $a \le t \le b$.
 Then for the inclusion $K : U(n) \subset Sp(2n,\mathbb{R})$:
 $$
 \begin{array}{l}

  Trace(  J \ (\frac{ d \ K(g(t))}{dt}) \cdot K(g(t))^{-1})
 = 2 \cdot \Sigma_{j=1}^n  \frac{ d\theta_j(t)}{dt}

 \end{array}
 $$
\end{thm}
This is to be compared to $(det_\mathbb{C} g(t))^2 = e^{2i \Sigma_{j=1}^n \ \theta_j(t)}$.
Applied to the cases two $G, \Pi: Sp(2n,\mathbb{R}) \rightarrow U(n)$,
this theorem evaluates the integrands
$$
( f \rightarrow  det_C(G[M(t)])^\star \ \frac{d\theta}{2 \pi} \ respectively \ \
( f \rightarrow  det_C(\Pi(M(t)))^\star \ \frac{d \theta}{2\pi}.
$$
arising from the unitary terms as $\frac{1}{2\pi}  Trace(  J \ (\frac{ d \ K(g(t))}{dt}) \cdot K(g(t))^{-1}) \
d\theta$ with $g(t)= G[M(t)]$ respectively $g(t)= \Pi(M(t))$.

The next theorem relates theorems \ref{thmGL1}, \ref{thmGL2}.
\begin{thm} \label{thm3}
For $f: [a,b] \rightarrow Sp(2n,\mathbb{R})$ smooth and decomposition as above,
$f(t) = K(\Pi(f(t))) \cdot e^{A(t)}$
$$ \begin{array}{l}
  Trace( J \ f(t)' \cdot f(t)^{-1}) =  Trace(  J \ K(U(f(t)))' \cdot K(U(f(t)))^{-1})\\ \hspace{1in}
  +  Trace( J \ A(t)' A(t)^{-1} \ (1/2) \ (Id[2n] - e^{-2A(t)}))
\end{array} $$
\end{thm}

The approach of CLM is related to the above by:
\begin{thm} \label{thm5}
Let $f: [a,b] \rightarrow Sp(2n,\mathbb{R})$ be a smooth symplectic path
with $f(t) = \left( \begin{array}{cc} a(t) & b(t) \\ c(t) & d(t) \end{array} \right)$.
Recall that \newline  $G[F(t)]= (a(t) \ + i \ c(t)) \cdot (a(t)^t a(t) + c(t)^t c(t))^{-1/2}$
is a unitary matrix and there is the decomposition
$
f(t) = K( G[f(t)]) \cdot L(t)
$
with
$$
L(t) = \left( \begin{array}{cc} D[M](t)^{+1/2} & 0 \\ 0 & D[M](t)^{-1/2} \end{array} \right)
\cdot \left( \begin{array}{cc} Id[n] & \alpha(t) \\ 0 & Id[n] \end{array} \right)
$$
where
$D[M](t)= (a(t)^t a(t) + c(t)^t c(t))$ and $\alpha(t) = D[M](t)^{-1} \ (a(t)^t b(t) + c(t)^t d(t))$
are both real symmetric.

Then
$$ \begin{array}{l}
Trace( J \ \frac{d f(t)}{dt} \ f(t)^{-1})
= Trace( J \ \frac{d K(G[f(t)])}{dt} \  K(G[f(t)])^{-1}\ )\\ \hspace{1in} - Trace( \ D[M](t) \ \frac{d \alpha(t)}{dt} \ ). \end{array}
$$
\end{thm}

\subsection{Details of Proofs for section \ref{sect9}:} \label{sect10}

\begin{lemma} \label{lemma1} Let $R,S : [a,b] \rightarrow Sp(2n,\mathbb{R})$ be smooth mappings,
then
$$ \begin{array}{l}
Trace( J \ (R(t)\cdot  S(t))' \cdot (R(t) \cdot S(t))^{-1})\\
 = Trace( J \ R(t)' \cdot R(t)^{-1} ) + Trace( (R(t)^{-1} \ J \ R(t))  \ S(t)' \cdot S(t)^{-1})
 \end{array}
 $$
 In particular, if $R(t)^t \cdot R(t) = Id[2n]$, then  $R(t)^{-1} \cdot J \cdot R(t)
 = R(t)^t \cdot J \cdot R(t) = J$ and
 $$ \begin{array}{l}
 Trace( J \ (R(t)\cdot  S(t))' \cdot (R(t) \cdot S(t))^{-1})\\
 = Trace( J \ R(t)' \cdot R(t)^{-1}) + Trace( J \ S(t)' \cdot S(t)^{-1}). \end{array}
 $$
 Here $R(t) \in Sp(2n,\mathbb{R}) \cap O(2n)$ is equivalent to $R(t) $ in the image of $K : U(n) \subset Sp(2n,\mathbb{R})$.
 \end{lemma}

 Proof lemma \ref{lemma1}: By $Trace(R \cdot S) = Trace(S \cdot R)$ and $Trace(k \ R + l \ S)=
 k \ Trace(R)+ l \ Trace(S)$:
 $$
 \begin{array}{l}
 Trace( J \ (R(t)\cdot  S(t))' \cdot (R(t) \cdot S(t))^{-1})\\
 = Trace( J \ ( R(t)' \cdot S(t) + R(t) \cdot S(t)') \cdot( S(t)^{-1} \cdot R(t)^{-1}) \\
 = Trace( J \ R(t)' \cdot R(t)^{-1}) + Trace( J \cdot R(t) \cdot S(t)'  \cdot( S(t)^{-1} \cdot R(t)^{-1}) \\
 = Trace(  J \ R(t)' \cdot R(t)^{-1}) + Trace( (R(t)^{-1} \ J \ R(t)) \cdot S(t)' \cdot B(t)^{-1})
 \end{array}
 $$
The lemma is proved.

\vspace{.3in}
Applying this lemma to theorem \ref{thm3} as $\Pi(f(t)) \in U(n)$
$$ \begin{array}{l}
Trace( J \ f(t)' \cdot f(t)^{-1}) \\ =
  Trace(  J \ K(\Pi(f(t)))' \cdot K(\Pi(f(t))^{-1}) + Trace( J \  (e^{A(t)})' \cdot e^{-A(t)} )
  \end{array}
  $$
  where and
remembering that $J \ A(t) = - A(t) \ J$:
$$\begin{array}{l}
  Trace( J \  (e^{A(t)})' \cdot e^{-A(t)} )\\
  = \Sigma_{n=1}^\infty \  (1/n!) \Sigma_{j=1}^n \  Trace( J,  \ replace \ j^{th} \ A(t) \ in \ A(t)^n\
  \ by \  A(t)') ) \cdot e^{-A(t)}\\
  and \ passing \ first \ j-1  \ A(t)'s  \ forwards \ past \ J \ then \ to \ back \ \ past \ e^{-A(t)}\\
 = \Sigma_{n=1}^\infty \  (1/n!) \Sigma_{j=1}^n \ (-1)^{j-1} \  Trace( J \ (A(t)' \cdot A(t)^{-1} )
  \cdot A(t)^n \cdot e^{-A(t)} )\\
  =  \Sigma_{n, odd} \   Trace(  J \ (A(t)' \cdot A(t)^{-1}) \ ( (1/n!) A(t)^n) \cdot e^{-A(t)} ) \\
  = Trace( J \ (A(t)' A(t)^{-1}) \ (1/2)(e^{A(t)} - e^{-A(t)}) \ e^{-A(t)})\\
= Trace( J \ A(t)' A(t)^{-1} \ (1/2) \ (Id[2n] -  e^{-2A(t)}))
  \end{array}
  $$

 Proof of theorem \ref{thm4}:

By application of lemma \ref{lemma1} via $U(t), V(t), Diag(t) \in U(n)$
and $(V(t)^{-1})'= - V(t)^{-1} \ V(t)' \ V(t)^{-1} $
$$
\begin{array}{l}
Trace( J \ K(U(t))' \ K(U(t))^{-1})\\ = Trace( J \ K( V(t) \ Diag(t) V(t)^{-1})' \ K( V(t) \ Diag(t) V(t)^{-1})^{-1})\\
= Trace( J \ K(V(t))' \ K(V(t))^{-1}) + Trace( J \ K(Diag(t))' \ K(Diag(t))^{-1})\\
\hspace{.6in} + Trace( J \ K(V(t)^{-1})' \ K(V(t)))\\
 = Trace(J \ K(V(t))' \ K(V(t))^{-1}) + Trace( J \ K(Diag(t))' \ K(Diag(t))^{-1})\\
\hspace{.6in}  - Trace( J \ K(V(t)^{-1} \ V(t)' )  )\\
  = Trace( J \ K(V(t))' \ K(V(t))^{-1}) + Trace( J \ K(Diag(t))' \ K(Diag(t))^{-1})\\
\hspace{.6in}  - Trace( J \ K(V(t))' \ K(V(t))^{-1} )\ by \ J, \ K(V(t)) \ commuting  \\
  = Trace( J \ K(Diag(t))' \ K(Diag(t))^{-1} =2 \cdot  \Sigma_{j=1}^n \ \frac{d \theta_j(t)}{dt}
  \end{array}
  $$
  To see the last equality, let $R(\theta(t)) = \left( \begin{array}{cc} cos(\theta(t)) & - sin(\theta(t)) \\ sin(\theta(t)) & cos(\theta(t)) \end{array} \right)$. Then
  $$
  \begin{array}{l}
    J \ \frac{R(\theta(t)}{ dt} \cdot R(\theta(t))^{-1}=  J \ \frac{R(\theta(t)}{ dt} \cdot R(-\theta(t))\\
   =  J \ \left( \begin{array}{cc} -sin(\theta(t)) & - cos(\theta(t)) \\ cos(\theta(t)) & -sin(\theta(t)) \end{array} \right) \left( \begin{array}{cc} cos(\theta(t)) & + sin(\theta(t)) \\- sin(\theta(t)) & cos(\theta(t)) \end{array} \right) \cdot \frac{d\theta}{dt}\\
   = \left( \begin{array}{cc} 0 & 1 \\ -1 & 0 \end{array} \right) \
    \left( \begin{array}{cc} 0 & -1 \\ +1 & 0 \end{array} \right) \ \frac{d\theta(t)}{dt}
    = \left( \begin{array}{cc} 1 & 0 \\ 0 & 1 \end{array} \right) \ \frac{d\theta(t)}{dt}

     \end{array}
   $$
Hence,
$$
Trace( J \ \frac{d R(\theta)}{dt} \ R(\theta(t))^{-1})
= 2 \ \frac{d\theta(t)}{dt}
$$

 \vspace{.1in}
 Proof of theorem \ref{thm5}:
 $
 f(t) = K(G[F(t)]) \cdot L(t)
 $
 with $L(t)$ specified above.

 Then, as before, by $G[f(t)]$ unitary so
 $$
 \begin{array}{l}
  \frac{ d f(t)}{dt} \ F(t)^{-1} \\
  = Trace( \ J \ ( \  \frac{dK(G[f(t)])}{dt} \ K(G[f(t)])^{-1})) + Trace(  J \ (\frac{d L(t)}{dt}  L(t)^{-1}))
  \end{array}
  $$

  It remains to evaluate the last term using $(V(t)^{-1})' = - V(t)^{-1} \ V(t)' \ V(t)^{-1}$.

  $$ \begin{array}{l}
  L(t) =  \left( \begin{array}{cc} D[M](t)^{+1/2} & 0 \\ 0 & D[M](t)^{-1/2} \end{array} \right)
\cdot \left( \begin{array}{cc} Id[n] & \alpha(t) \\ 0 & Id[n] \end{array} \right)\\

\hspace{1in} = \left( \begin{array}{cc} D[M](t)^{+1/2} & D[M](t)^{1/2} \alpha(t) \\ 0 & D[M](t)^{-1/2} \end{array} \right)\\

So , \
   L(t)^{-1} =\left( \begin{array}{cc} Id[n] & -\alpha(t) \\ 0 & Id[n] \end{array} \right) \cdot  \left( \begin{array}{cc} D[M](t)^{-1/2} & 0 \\ 0 & D[M](t)^{+1/2} \end{array} \right)

 \\ \hspace{1in} =  \left( \begin{array}{cc} D[M](t)^{-1/2} & -\alpha(t) \ D[M](t)^{1/2} \\ 0 &  D[M](t)^{1/2} \end{array} \right), \  \\

   and \  \ J \  \frac{d L(t)}{dt} \ L(t)^{-1} \\
 = J \ \left( \begin{array}{cc} \frac{ d D[M](t)^{+1/2} }{dt} &  \frac{d (D[M](t)^{1/2} \alpha(t))}{dt} \\ 0 &- D[M](t)^{-1/2} \
 \frac{d D[M](t)^{1/2}}{dt} \
 D[M](t)^{-1/2} \end{array} \right)\\ \hspace{1in} \cdot  \left( \begin{array}{cc} D[M](t)^{-1/2} & - \alpha(t) D[M](t)^{1/2} \\ 0 & D[M](t)^{+1/2} \end{array} \right) \\
J \  \left( \begin{array}{cc} \frac{D[M](t)^{1/2}}{dt} \ D[M](t)^{-1/2} &
Q  \\ 0 & -D[M](t)^{-1/2} \ \frac{d D[M](t)^{1/2}}{dt} \end{array} \right)\\
\end{array} $$

 $$
 \begin{array}{l}
 with \ Q =
- \frac{d D[M](t)^{1/2}}{dt} \alpha(t) D[M](t)^{1/2}  + \frac{d ( D[M](t)^{1/2} \alpha(t) )}{dt} D[M](t)^{1/2},  \\
 \hspace{1in} =  D[M](t)^{1/2} \ \frac{d  \alpha(t) }{dt} D[M](t)^{1/2}\\
 So \ Trace(J \  \frac{d L(t)}{dt} \ L(t)^{-1}) \\
=Trace  \left( \begin{array}{cc} 0 & 1 \\ -1 & 0 \end{array} \right) \  \left( \begin{array}{cc} \frac{D[M](t)^{1/2}}{dt} \ D[M](t)^{-1/2} &

  Q \\ 0 & -D[M](t)^{-1/2} \ \frac{d D[M](t)^{1/2}}{dt} \end{array} \right)\\
  =
  Trace(\left( \begin{array}{cc} 0 & -D[M](t)^{-1/2} \ \frac{d D[M](t)^{1/2}}{dt} \\ - \frac{D[M](t)^{1/2}}{dt} \ D[M](t)^{-1/2} &

 - Q  \end{array} \right))\\
 = - Trace(Q) = - Trace( D[M](t)^{1/2} \ \frac{d  \alpha(t) }{dt} D[M](t)^{1/2})
 = - Trace( D[M](t) \ \frac{d \alpha(t)}{dt})
 \end{array}
 $$

\vspace{-0.1in}

\frenchspacing

\bibliographystyle{plain}

\vspace{.3in}

Courant Institute of Mathematics, New York University, \ New York, NY 10012

email addresses: cappell@cims.nyu.edu,  em1613@nyu.edu

\end{document}